\documentclass[a4paper,11pt,DIV=11,
abstract=on,
]{scrartcl}
\usepackage{mathtools}
\mathtoolsset{showonlyrefs}
\usepackage{amsmath, amsthm, amssymb}
\usepackage{fontenc}
\usepackage[utf8]{inputenc}
\usepackage[english]{babel}
\usepackage{graphicx}
\usepackage{subcaption,xargs}
\usepackage{algorithm, booktabs}
\usepackage[noend]{algpseudocode}
\usepackage{enumitem,listings,microtype}
\usepackage{environ,ifthen,xcolor}
\usepackage{multirow}
\usepackage[
  backend=bibtex,
  style=numeric-comp,
  firstinits=true,
  natbib=true,
  maxcitenames=5,
  sortlocale=en_EN,
  url=true, 
  doi=true,
  eprint=true
]{biblatex}
\usepackage{hyperref}

\DeclareCaptionLabelSeparator{periodspace}{.\ }
\setlist[enumerate,1]{label={\alph*)}}
\setlist[enumerate,2]{label={\roman*)}}

\setkomafont{sectioning}{\rmfamily\bfseries}
\setkomafont{title}{\rmfamily}
\newcommand{\e}{{\,\mathrm{e}}}
\newcommand{\im}{{\mathrm{i}}}
\newcommand{\p}{{\mathrm{p}}}

\newcommand{\dx}{\,\mathrm{d}}
\makeatletter
\newcommandx{\abs}[2][1=\@empty]{#1\lvert #2 #1\rvert}
\newcommandx{\norm}[3][1=\@empty,3=\@empty]{#1\lVert #2 #1\rVert_{#3}}
\makeatother
\newcommand{\tT}{\mathrm{T}}

\newcommand{\vect}[1]{\mathbf{#1}}
\newcommand{\mat}[1]{\mathbf{#1}}
\newcommand{\Stiffness}{\mathcal{C}}
\newcommand{\eff}{\mathrm{eff}}

\newcommand{\Z}{\mathbb{Z}}
\newcommand{\R}{\mathbb{R}}

\newcommand{\T}{\mathbb{T}}
\newcommand{\Wiener}{A}
\newcommand{\tensorProd}{:}
\newcommand{\Def}[1]{\emph{#1}}
\newcommand{\Curlfree}{\mathcal{E}}

\DeclareMathOperator{\Div}{div}
\DeclareMathOperator{\Id}{Id}
\DeclareMathOperator{\Fourier}{\mathcal F}
\DeclareMathOperator{\generatingSet}{\mathcal G}
\DeclareMathOperator{\Grad}{\nabla}
\DeclareMathOperator{\Pattern}{\mathcal P}
\DeclareMathOperator{\GradSym}{\Grad_{\mathrm{Sym}}}
\DeclareMathOperator{\spanOp}{span}
\DeclareMathOperator{\sinc}{sinc}
\DeclareMathOperator{\Sym}{Sym}
\DeclareMathOperator{\SSym}{SSym}
\DeclareMathOperator{\Translate}{\mathcal T}
\DeclareMathOperator{\diag}{diag}

\DeclareMathOperator{\Fundamental}{I}
\hypersetup{pdfauthor={Ronny Bergmann, Dennis Merkert},
	pdftitle={FFT-based homogenization on periodic anisotropic translation invariant
	spaces},
	pdfsubject={preprint},%
	pdfcreator = {pdflatex and TextMate},%
	pdfkeywords={},
	plainpages=false, pdfstartview=FitH, pdfview=FitH, pdfpagemode=UseOutlines,%
	bookmarksnumbered=true,bookmarksopen=false,bookmarksopenlevel=0,%
	colorlinks=true,linkcolor=black,citecolor=black,urlcolor=black%
}

\newtheorem{theorem}{Theorem}[section]
\newtheorem{definition}[theorem]{Definition}
\newtheorem{proposition}[theorem]{Proposition}
\newtheorem{lemma}[theorem]{Lemma}
\newtheorem{remark}[theorem]{Remark}

\title{FFT-based homogenization on periodic anisotropic translation invariant
	spaces}

\author{%
Ronny Bergmann\thanks{%
Department of Mathematics,
University of Kaiserslautern,
Postfach 3049, D-67653 Kaiserslautern, Germany\newline$\{$bergmann, dmerkert$\}$@mathematik.uni-kl.de.
}
\and
Dennis Merkert\footnotemark[1]}
\date{January 17, 2017}
\begin{document}
\maketitle
\begin{abstract}
	\noindent\small
	In this paper we derive a discretisation of the equation of quasi-static
	elasticity in homogenization in form of a variational formulation and
	the so-called Lippmann-Schwinger equation, in anisotropic spaces of
	translates of periodic functions. We unify and extend the truncated
	Fourier series approach, the constant finite element ansatz and the
	anisotropic lattice derivation. The resulting formulation of the
	Lippmann-Schwinger equation in anisotropic translation invariant spaces
	unifies and analyses for the first time both the Fourier methods and
	finite element approaches in a common mathematical framework. We further
	define and characterize the resulting periodised Green operator. This
	operator coincides in case of a Dirichlet kernel corresponding to a
	diagonal matrix with the operator derived for the Galerkin projection
	stemming from the truncated Fourier series approach and to the
	anisotropic lattice derivation for all other Dirichlet kernels.
	Additionally, we proof the boundedness of the periodised Green operator.
	The operator further constitutes a projection if and only if the space
	of translates is generated by a Dirichlet kernel. Numerical examples for
	both the de la Vall\'ee Poussin means and Box splines illustrate the
	flexibility of this framework.
\end{abstract}
%
%
\section{Introduction}\label{sec:Introduction}
Many modern tools and products use composites, i.e. mixtures of two or more
materials with distinct elastic properties to obtain certain flexible behaviour,
dampening effects, or longevity. Homogenization aims to simplify simulations by
replacing the microscopically composed material by a homogeneous one which
behaves the same on the macroscopic scale. Mathematically one assumes is a
periodic microstructure, i.e.~a structure that can be represented by a certain
unit cell with periodic boundary conditions.

For the simulation of such elastic composite structures Moulinec and
Suquet~\cite{MoulinecSuquet1994,MoulinecSuquet1994} derive an algorithm based
on the fast Fourier transform. This algorithm, called the Basic Scheme, inspired
many similar numerical approaches based on using discretised differential
operators~\cite{Willot2015,Schneider2015,Schneider2016} and extensions to porous
media~\cite{Michel2000,Schneider2016}. Information on sub-structures of the
geometry is incorporated into the solution
method in~\cite{KMS:2015Homogenization}.
The Basic Scheme is generalized to problems of higher order, i.e.
derivatives of strain and stiffness~\cite{Tran2012}, and the solution of the
arising linear system by Krylov subspace methods is
analysed in~\cite{ZVNM:2010NumHom}.

Vond\v{r}ejc et.al.~\cite{Vondrejc2014} show that the method of Moulinec and
Suquet can also be understood as a Galerkin projection using truncated Fourier
series. This idea is generalized in~\cite{BergmannMerkert2016} to anisotropic
lattices thus allowing to take directional information on the geometrical
structure or the orientation of interfaces between materials into account.
Brisard and Dormieux~\cite{BrisardDormieux2010Framework,BrisardDormieux2012} use
constant finite elements to arrive at the Basic Scheme with a modified linear
operator, based on an energy based formulation.

In this paper we unify and extend the approaches of Vond\v{r}ejc et.al., Brisard
and Dormieux, and the anisotropic lattice ansatz obtaining a discretisation of
the equations for quasi-static elasticity in homogenization in anisotropic
spaces of periodic translates. Vond\v{r}ejc et.al. show that a variational
equation and a formulation with the strain as a fixed-point, the so-called
Lippmann-Schwinger equation, are equivalent by means of a projection operator
derived from the Green operator. This paper introduces a periodised Green
operator on the space of translates for the Lippmann-Schwinger equation.
Furthermore we classify the properties of this operator in case of spaces of
translates and prove that it induces a projection operator if and only if the
space of translates is generated by the Dirichlet kernel. Hence the
(anisotropic) truncated Fourier series emerges as a case with special properties
of the general setting introduced here. This introduces further insight into the
equivalence of the variational and the Lippmann-Schwinger equation formulation
of Vond\v{r}ejc et.al. The mathematical framework this paper introduces unifies
and analyses the
approaches of Fourier methods and finite elements for the first time. Especially
using translates of Box splines~\cite{deBoorHoelligRiemenscheider1993BoxSplines}
as ansatz functions incorporates the constant finite elements of Brisard and
Dormieux allowing also for anisotropic finite elements of arbitrary smoothness.
A different approach to solve the equation using linear finite elements with full
quadrature is shown in~\cite{Schneider2016} and is based on replacing the
continuous differential operator by a discrete one.

Spaces of translates can for example be generated by de la Vall\'ee Poussin
means which provide a generalization of the Dirichlet and Fej\'er kernel. They
combine a finite support in frequency domain with good localization in
space~\cite{GohGoodman2004}. These functions introduce a trade off between
damping of the Gibbs phenomenon and reproduction of multivariate trigonometric
monomials. They allow for better predictions of the elastic macroscopic
properties of the composite material and result in smoother --- and thus better
--- solutions. Further, different Box splines and their influence on the
solution are demonstrated.

The remainder of the paper is structured as follows. After reviewing important
properties of anisotropic spaces of translates in
Section~\ref{sec:Preliminaries} the partial differential equation of
quasi-static elasticity in homogenization is introduced in
Section~\ref{sec:homogenization}. Then, the periodised Green operator on spaces
of translates is introduced. This operator is subsequently analysed regarding
its equivalence to a projection operator and then used to discretise the
Lippmann-Schwinger equation. Numerical examples are then provided in
Section~\ref{sec:numerics} and a conclusion is drawn in
Section~\ref{sec:summary}.
%
%
\section{Preliminaries}\label{sec:Preliminaries}

Throughout this paper we will employ the following notation: the symbols
\(a\in\mathbb C\), \(\vect{a}\in\mathbb C^d\) and \(\mat{A}\in\mathbb C^{d\times
d}\) denote scalars, vectors, and matrices, respectively. The only exception
from this are~\(f,g,h\) which are reserved for functions. We denote the inner
product of two vectors by \(\vect{a}^\tT\vect{b} \coloneqq\sum_i a_ib_i\) and
reserve the symbol \(\langle\cdot,\cdot\rangle\) for inner products of two
functions or two generalized sequences, respectively. For a complex number \(a =
b+\im c\), \(b,c\in\mathbb R\), we denote the complex conjugate by
\(\overline{a} \coloneqq b-\im c\). Constants like Euler's number~\(\e\) or the
imaginary unit \(\im\), i.e.\ \(\im^2 = -1\), are set upright.

Usually, we are concerned with \(d\)-dimensional data, where \(d=2,3\), but the
theory is written in arbitrary dimensions. Sets are denoted by capital
case calligraphic letters, e.g.~\(\mathcal P\) or~\(\mathcal G\) and the same
for the Fourier transform~\(\mathcal F\) which all might depend on parameters
given in round brackets. We denote second-order tensors by small Greek letters
as~\(\lambda,\varepsilon\) with entries \(\lambda_{ij}\) are indexed again by
scalars \(i,j\) and similarly we denote fourth-order tensors by capital
calligraphic letters, where \(\mathcal C\) is the most prominent one.

\subsection{Arbitrary patterns and the Fourier transform}
\label{subsec:preliminariesFourier}
The space of functions we are concerned with is the Hilbert
space~$L^2(\T^d)$ of (equivalence classes of) square integrable functions
on the \(d\)-dimensional torus~$\T \cong [-\pi,\pi)^d$ with inner product
\begin{equation*}
	\langle f,g \rangle 
	=\frac{1}{(2\pi)^d}
		\int_{\T^d}
			f(\vect{x})\overline{g(\vect{x})}
		\dx\vect{x},
		\qquad
		f,g \in L^2(\T^d)
	\text{.}
\end{equation*}
In several cases, the functions of interest are tensor-valued.
For these functions, we take the tensor product of the Hilbert space,
e.g.~\(L^2(\mathbb T^d)^{n\times n}\) for the space of
functions~\(f\colon \mathbb T^d \to \mathbb C^{n\times n}\) that have values
being~\(n\times n\)-dimensional matrices. The following
preliminaries can be generalized to these tensor product spaces by performing
the operations element wise. We restrict the following of this
subsection therefore to the case of \(L^2(\mathbb T^d)\).

Every function $f \in L^2(\T^d)$ can be written in its Fourier series
representation
\begin{equation}\label{eq:fourier-series}
	f(\vect{x})
	= \sum_{\vect{k} \in \mathbb Z^d} c_{\vect{k}}(f)\e^{\im\vect{k}^\tT\vect{x}},
\end{equation}
introducing the multivariate Fourier
coefficients~\(c_{\vect{k}}(f) = \langle f,\e^{\im\vect{k}^\tT\circ}\rangle\),
\(\vect{k}\in\mathbb Z^d\). The equality in \eqref{eq:fourier-series} is meant in  \(L^2(\T^d)\) sense.
We denote by $\vect{c}(f)
	= \bigl\{c_{\vect{k}}(f)\bigr\}_{\vect{k}\in\mathbb Z^d}
	\in \ell^2(\mathbb Z^d)$ generalized sequences which form a Hilbert
space with the inner product
\begin{equation*}
	\langle\vect{c},\vect{d}\rangle
		= \sum_{\vect{k} \in \mathbb Z^d}c_{\vect{k}}\overline{d_{\vect{k}}},
		\qquad \vect{c},\vect{d}\in \ell^2(\mathbb Z^d)\text{.}
\end{equation*}
The Parseval equation reads
\begin{equation}
		\langle f, g \rangle
		= \langle \vect{c}(f),\vect{c}(g) \rangle
		= \sum_{\vect{k} \in \mathbb Z^d} c_{\vect{k}}(f) \overline{c_{\vect{k}}(g)}\text{.}
		\label{eq:parseval}
\end{equation}

\paragraph{The pattern and the generating set.}
For any regular matrix $\mat{M} \in \mathbb Z^{d\times d}$ we define the
congruence relation for $\vect{h},\vect{k} \in \mathbb Z^d$ with respect
to~$\mat{M}$ by
\begin{equation*}
		\vect{h} \equiv \vect{k} \bmod \mat{M}
		\Leftrightarrow \exists\,\vect{z} \in \mathbb Z^d\colon \vect{k} = \vect{h} + \mat{M}\vect{z}\text{.}
\end{equation*}
We define the lattice
\[
	\Lambda(\mat{M}) \coloneqq \mat{M}^{-1}\mathbb Z^d
	= \{\vect{y}\in\mathbb R^d : \mat{M}\vect{y} \in \mathbb Z^d\},
\]
and the pattern \(\Pattern(\mat{M})\) as any set of congruence representants of
the lattice with respect to \(\bmod\ \vect{1}\), e.g.\ $ \Lambda(\mat{M})\cap[0,1)^d$
or~$\Lambda(\mat{M})\cap\bigl[-\tfrac{1}{2},\tfrac{1}{2}\bigr)^d$. For the rest
of the paper we will refer to the set of congruence class representants in the 
symmetric unit cube \(\bigl[-\tfrac{1}{2},\tfrac{1}{2}\bigr)^d\). The generating
set \(\generatingSet(\mat{M})\) is defined by \(\generatingSet(\mat{M}) 
\coloneqq \mat{M}\Pattern(\mat{M})\) for any pattern \(\Pattern(\mat{M})\). For
both, the number of elements is given by \(
	\abs{\Pattern(\mat{M})}
	=\abs{\generatingSet(\mat{M})}
	=\abs{\det{\mat{M}}}
	\eqqcolon m,
\) which follows directly from~\cite[Lemma II.7]{deBoorHoelligRiemenscheider1993BoxSplines}.

For a regular integer matrix \(\mat{M}\in\mathbb Z^{d\times d}\) and an
absolutely summable generalized sequence \(\vect{a} =
  \{a_{\vect{k}}\}_{\vect{k}\in\mathbb Z^d}\) we further define the
\emph{bracket sum}
\begin{equation}\label{eq:bracketsum}
\bigl[\vect{a}\bigr]_{\vect{k}}^{\mat{M}}
\coloneqq
\sum_{\vect{z}\in\mathbb Z^d} a_{\vect{k}+\mat{M}^\tT\vect{z}},
\qquad \vect{k}\in\mathbb Z^d.
\end{equation}
The bracket sum is periodic with respect to \(\mat{M}^\tT\),
i.e.,~\(\bigl[\vect{a}\bigr]_{\vect{k}}^{\mat{M}} = \bigl[\vect{a}\bigr]_{\vect{k}+\mat{M}^\tT\vect{z}}^{\mat{M}}\) holds for any \(\vect{k},\vect{z}\in\mathbb Z^d\).

\paragraph{A fast Fourier transform on patterns.}
The discrete Fourier transform on the pattern $\Pattern(\mat{M})$ is defined~\cite{ChuiLi:1994} by
\begin{equation}\label{eq:Fouriermatrix}
	\mathcal F(\mat{M})
	\coloneqq
	\frac{1}{\sqrt{m}}
	\Bigl(
		\e^{- 2\pi \im \vect{h}^\tT\vect{y}}
	\Bigr)_{%
		\vect{h} \in \generatingSet(\mat{M}^\tT),\,%
		\vect{y} \in \Pattern(\mat{M})%
},
\end{equation}
where $\vect{h}\in\generatingSet(\mat{M}^\tT)$ indicate the rows and
$\vect{y} \in \Pattern(\mat{M})$ indicate the columns of the Fourier
matrix~\(\mathcal F(\mat{M})\).
The discrete Fourier transform
on $\Pattern(\mat{M})$ is defined for a
vector~$\vect{a} = (a_{\vect{y}})_{\vect{y}\in\Pattern(\mat{M})}\in\mathbb C^m$
arranged in the same ordering as the columns in~\eqref{eq:Fouriermatrix} by
\begin{equation}\label{eq:FourierTransform}
	\vect{\hat a} = (\hat a_{\vect{h}})_{\vect{h}\in\generatingSet(\mat{M}^\tT)}
	= \mathcal F(\mat{M})\vect{a},
\end{equation}
where the resulting vector $\vect{\hat a}$ is ordered as the columns
of~$\mathcal F(\mat{M})$ in~\eqref{eq:Fouriermatrix}. Its
implementation yields complexity of~\(\mathcal O(m\log m)\) similar to the
classical Fourier transform, when the ordering is fixed as described
in~\cite[Theorem~2]{Bergmann2013FFT}. Note that for the so-called rank-1-lattices, the Fourier transform on the pattern even reduces
to a one-dimensional FFT for patterns in arbitrary
dimensions~\cite{KaemmererPottsVolkmer2015b}.
\subsection{Translation invariant spaces of periodic functions}%
\label{subsec:preliminariesdlVP}
\paragraph{Spaces of translates and interpolation.}
	A space of functions \(V\subset L^2(\T^d)\) is called \(\mat{M}\)-invariant,
	if for all \(\vect{y}\in\Pattern(\mat{M})\) and all functions \(f\in V\) the
	translates \(\Translate({\vect{y}})f\coloneqq f(\cdot-2\pi\vect{y})\in V\). Especially the space
	\[
		V_{\mat{M}}^{f} \coloneqq \spanOp\bigl\{\Translate(\vect{y})f\,:\,
		\vect{y}\in\Pattern(\mat{M})\bigr\}
	\]
	of translates of \( f \) is \( \mat{M} \)-invariant.
	A function \(g\in V_{\mat{M}}^f \) is of the form
	\begin{equation*}
		g=\displaystyle\sum_{\vect{y}\in\Pattern(\mat{M})}
			a_{\vect{y}}\Translate(\vect{y})f
	\end{equation*}
	For \( f\in L^2(\T^d) \) an easy calculation on the Fourier
	coefficients using the unique decomposition of \( \vect{k}\in\mathbb Z^d \)
	into \( \vect{k} = \vect{h}+\mat{M}^\tT\vect{z} \),
	\( \vect{h}\in\generatingSet(\mat{M}^\tT), \vect{z}\in\mathbb Z^d \),
	yields, that \( g\in V_{\mat{M}}^f \) holds if and only if~\cite[Theorem 3.3]{LangemannPrestin2010WaveletAnalysis}
  \begin{equation}\label{eq:inTranslates:ck}
		c_{\vect{h}+\mat{M}^\tT\vect{z}}(g) 
		= \hat a_{\vect{h}}c_{\vect{h}+\mat{M}^\tT\vect{z}}(f)
		\quad\text{for all }
		\vect{h}\in\generatingSet(\mat{M}^\tT), \vect{z}\in\mathbb Z^d\text{,}
	\end{equation}
	holds, where
	\(
	\vect{\hat a}
	= \bigl(\hat a_{\vect{h}}\bigr)_{\vect{h}\in\generatingSet(\mathbf{M}^\tT)}
	= \sqrt{m}\mathcal F(\mat{M})\vect{a}
	\)
	denotes the discrete Fourier transform
	of~\(
		\vect{a}
		= \bigl( a_{\vect{y}}\bigr)_{\vect{y}\in\Pattern(\mathbf{M})}
			\in\mathbb C^m
		      \), see~\cite{LangemannPrestin2010WaveletAnalysis}.
	Using the space of trigonometric polynomials on the
	generating set \( \generatingSet(\mat{M}^\tT) \), which is denoted by
	\[
	\Translate_{\mat{M}} \coloneqq\Bigl\{f\,:\,
	f=
	\sum_{\vect{h}\in\generatingSet(\mat{M}^\tT)}a_{\vect{h}}\e^{\im\vect{h}^\tT\circ},\ a_{\vect{h}}\in\mathbb C
	\Bigr\}\text{,}
	\]
	we define for a function \( f\in L^2(\T^d) \) the Fourier partial
	sum~\( \operatorname{S}_{\mat{M}}f\in\Translate_{\mat{M}} \) by
	\[
		\operatorname{S}_{\mat{M}}f
		\coloneqq
		\sum_{\vect{h}\in\generatingSet(\mat{M}^\tT)}
		c_{\vect{h}}(f)\e^{\im\vect{h}^\tT\circ}.
	\]
	
	The discrete Fourier coefficients \(c_{\vect{k}}^{\mat{M}}(f)\) of a function
  \(f\) that is evaluated pointwise on the pattern \(\Pattern(\mat{M})\) are
  defined by
	\begin{equation*}
		c_{\vect{h}}^{\mat{M}}(f)
		\coloneqq
		\frac{1}{m}
			\sum_{\vect{y}\in\Pattern(\mat{M})}
				f(2\pi\vect{y})\e^{-2\pi\im\vect{h}^\tT\vect{y}}
		,\quad\vect{h}\in\generatingSet(\mat{M}^\tT)
		\text{.}
	\end{equation*}
	The discrete Fourier coefficients \(c_{\vect{h}}^{\mat{M}}\) are related to
	the Fourier coefficients for a function~\( f\in A(\T^d) \), where
	\(A(\T^d)\) denotes the Wiener Algebra, i.e. the space of functions
	with absolutely convergent Fourier series. This relation is given in the 
	following Lemma, also known as the aliasing formula, see
	e.g.~\cite[Lemma 2]{BergmannPrestin2014Interpolation}.
	\begin{lemma}\label{lem:Aliasing}
		Let \( f\in A(\T^d) \) and the regular
		matrix~\(\mat{M}\in\mathbb Z^{d\times d}\) be given.
		Then the discrete Fourier coefficients~\( c_{\vect{h}}^{\mat{M}}(f) \)
		are given by
		\begin{equation*}
			c_{\vect{h}}^{\mat{M}}(f) 
			= 
			\sum_{\vect{z}\in\mathbb Z^d} c_{\vect{h}+\mat{M}^\tT\vect{z}}(f)
			= \bigl[\vect{c}(f)\bigr]_{\vect{h}}^{\mat{M}}
			,\quad \vect{h}\in\generatingSet(\mat{M)^\tT}\text{.}
		\end{equation*}
	\end{lemma}
	When looking at the space~\(V_{\mat{M}}^f\) of translates, the following
	definition is crucial in order to approximate a function \(g\) by using
	these translates.
	\begin{definition}
		\label{def:IP}
		Let \( \mat{M}\in\mathbb Z^{d\times d} \) be a regular matrix.
		A function~\( \operatorname{I}_\mat{M} \in V_\mat{M}^\varphi\) is
		called fundamental interpolant or Lagrange function
		of~\( V_\mat{M}^\varphi \) if
		\begin{equation*}
			\operatorname{I}_\mat{M}(2\pi\vect{y})
			\coloneqq \delta_{\vect{0},\vect{y}}^{{\mathbf E}_d},
			\quad \vect{y}\in\Pattern(\mat{M}),\quad\text{where } 
			\delta_{\vect{x},\vect{y}}^{{\mathbf M}} \coloneqq
			\begin{cases}
			  1, &\text{ if } \vect{y}\equiv \vect{x}\bmod{\mathbf M},\\
			  0, &\text{ else.}
			\end{cases}
		\end{equation*}
	\end{definition}
	The following lemma characterizes the existence of such a fundamental
	interpolant in a space \(V_\mat{M}^f\) of translates and collects some properties of
	the translates themselves, see~\cite[Lemma 1.23]{Bergmann2013Thesis} and~\cite[Lemma 2]{BergmannPrestin2014Interpolation}.
	\begin{lemma}\label{lem:FI:Props}
		Given a regular matrix \(\mat{M}\in\mathbb Z^{d\times d}\) and a
		function \(f\in A(\T^d)\), then the following holds.
		\begin{enumerate}
			\item \label{lem:FI:Prosp:Existence}The fundamental
			interpolant~\(\operatorname{I}_\mat{M}\in V_\mat{M}^f\) exists if
			and only if
			\begin{equation*}
				\bigl[\vect{c}(f)\bigr]_{\vect{h}}^{\mat{M}}
				\neq 0\quad\text{ for all }\vect{h}\in\generatingSet(\mat{M}^\tT).
			\end{equation*}
			If the fundamental
			interpolant~\( \operatorname{I}_\mat{M}\in V_\mat{M}^f \) exists,
			it is uniquely determined.
			\item\label{lem:FI:Props:LI} The set of
			translates~\(\bigl\{\Translate(\vect{y})f:\ \vect{y}\in\Pattern(\mat{M})\bigr\}\)
			is linear independent if and only if
			\[
				\sum_{\vect{z}\in\Z^d}
				\abs{c_{\vect{h}+\mat{M}^\tT\vect{z}}(f)}^2 >0
				\quad\text{holds for all }\vect{h}\in\generatingSet(\mat{M}^\tT).
			\]
			\item\label{lem:FI:Props:ONB} The set of
			translates~\(\bigl\{\Translate(\vect{y})f:\ \vect{y}\in\Pattern(\mat{M})\bigr\}\)
			is an orthonormal basis of \(V_{\mat{M}}^{f}\) if and only if
			\[
				\sum_{\vect{z}\in\Z^d}
				\abs{c_{\vect{h}+\mat{M}^\tT\vect{z}}(f)}^2 =
				\frac{1}{m}
				\quad\text{holds for all }\vect{h}\in\generatingSet(\mat{M}^\tT).
			\]
			\item\label{lem:FI:Props:IP} Given a
			function~\(\tilde g\in A(\T^d)\) we can obtain a
			function~\( g\in V_{\mat{M}}^{f}\) fulfilling
			\[
				\tilde g(2\pi\vect{y}) = g(2\pi\vect{y}),\quad \vect{y}\in\Pattern(\mat{M}),
			\]
			provided that the fundamental interpolant exists (which also
			implies linear independence of the translates on \(f\)) as
			\begin{equation*}
			  \hat a_{\vect{h}}
			  = \frac{%
			    \bigl[\vect{c}(\tilde g)\bigr]_{\vect{h}}^{\mat{M}}%
			    }{%
			    \bigl[\vect{c}(f)\bigr]_{\vect{h}}^{\mat{M}}%
			  },\qquad \vect{h}\in\generatingSet(\mat{M}^\tT),
			\end{equation*}
			where the coefficients~\(\hat a_{\vect{h}}\)
			yield~\(g\) in Fourier coefficients by~\eqref{eq:inTranslates:ck}.
		\end{enumerate}
	\end{lemma}
	By using Lemma~\ref{lem:FI:Props} changing from sampling values, i.e.~the
	coefficients on the pattern~\(\Pattern(\mat{M})\) of the fundamental
	interpolant, to coefficients with respect to~\(f\) in the corresponding
	space~\(V_{\mat{M}}^f\) of translates can be done by using the Fourier
	transform~\eqref{eq:FourierTransform} and the Fourier
	coefficients~\(c_{\vect{k}}(f)\) of~\(f\).

For the remainder of this paper, two special spaces of translates are of
interest, the periodised Box splines and the de la Vall\'ee Poussin means.

\paragraph{Periodised Box splines.}
\label{par:periodized_box_spline}
Let \(\mat{\Xi} = (\vect{\xi}_1,\ldots,\vect{\xi}_s)\in\mathbb R^{d\times s}\)
denote a set of column vectors \(\vect{\xi}_i\), \(i=1,\ldots,s\), where we
assume that these vectors span the \(\mathbb R^d\), i.e. especially we have
\(s\geq d\). Then the \emph{centred Box spline \(B^c_{\mat{\Xi}}\)}
can be defined via its Fourier transform as
\[
  \widehat B^c_{\mat{\Xi}}(\vect{y}) = \prod_{\vect{\xi}\in\mat{\Xi}} \sinc\bigl(\frac{1}{2}\vect{\xi}^\tT\vect{y}\bigr),
  \qquad
  \sinc(t) \coloneqq \frac{\sin(t)}{t},
\]
cf.~\cite[p. 11]{deBoorHoelligRiemenscheider1993BoxSplines}. A Box spline has compact support. For a function \(g\colon\mathbb R^d\to\mathbb C\) we can introduce its periodisation
\[
  g_{\mathrm{p}}(\vect{x}) \coloneqq \sum_{\vect{z}\in\mathbb Z^d} g\bigl(\frac{\vect{x}}{2\pi} - \vect{z}\bigr).
\]
Its Fourier coefficients can be directly computed from the continuous Fourier transform \(\hat g\) of \(g\), cf. e.g.~\cite[p.~41]{BergmannPrestin2014Interpolation}, as
\[
  c_{\vect{k}}(g_{\mathrm{p}}) = \hat g(2\pi\vect{k}),\quad \vect{k} \in \Z^d
\]
We combine these two to introduce the~\emph{periodised Box Spline} \(B_{\mat{\Xi}}\colon\mathbb T^d\to\mathbb R\) via its Fourier coefficients as
\[
  c_{\vect{k}}(B_{\mat{\Xi}}) \coloneqq 
  \widehat B^c_{\mat{\Xi}}(2\pi\vect{k}) = \prod_{\xi\in\Xi} \sinc\bigl(\pi\xi^\tT \vect{k}\bigr),\qquad \vect{k}\in\mathbb Z^d.
\]
Finally, we obtain by scaling the \emph{periodised pattern Box Spline} \(f_{\mat{M},\mat{\Xi}}\)
\[
  f_{\mat{M},\mat{\Xi}}(\vect{x}) \coloneqq B_{\mat{\Xi}}(\mat{M}^{-1}\vect{x}).
\]
Note that its translates might not be linearly independent for an arbitrary set
of vectors in~\(\mat{\Xi}\), see also~\cite{Poeplau1995}.
However, by~\cite{BergmannPrestin2014Interpolation}, see
also~\cite[Sect.~4]{deBoorHoelligRiemenschneider1985Interpolation}, the matrices
\(\mat{\Xi}\in\mathbb R^{d\times (p+q+r)}\) of the form
\(\xi_1=\ldots,\xi_p=(1,0)^\tT\), \(\xi_{p+1}=\ldots=\xi_{p+q}=(0,1)^\tT\), and
\(\xi_{p+q+1}=\ldots=\xi_{p+q+r}=(1,1)^\tT\), where at least two of the values
\(p,q,r\) are larger than \(0\), induce a periodised pattern Box spline
$f_{\mat{M},\Xi}$ with linear independent translates.

\paragraph{De la Vall\'ee Poussin means.} A special case of \(\mat{M}\)-invariant
spaces are the ones defined via de la Vall\'ee Poussin
means, following the construction of~\cite{BergmannPrestin2014dlVP}.
We call a function \(g\colon\mathbb R^d\to\mathbb R\) \emph{admissible} if the
function fulfils
\begin{enumerate}
	\item \(g(\vect{x}) \geq 0\) for all \(\vect{x}\in\R^d\),
	\item \(g(\vect{x}) > 0\) for \(\vect{x}\in\bigl[-\tfrac{1}{2},\frac{1}{2})^d\),
	\item \(\displaystyle\sum_{\vect{z}\in\Z ^d} g(\vect{x}+\vect{z}) = 1\)
	for all \(\vect{x}\in\R^d\).
\end{enumerate}
This can be for example the Box
splines of the form
\begin{equation*}
	g_{\vect{\alpha}}(\vect{x})
	\coloneqq B_{\mat{\Xi}}(\vect{x})
	,\quad\mat{\Xi} = \begin{pmatrix}
		\diag(\vect{\alpha}) \mat{I}_d
	\end{pmatrix}\in\mathbb R^{d\times 2d},\ \vect{\alpha}\in[0,1]^d
\end{equation*}
where \(\mat{I}_d\) is the \(d\)-dimensional unit matrix. 
In the following we define the de la Vall\'ee Poussin means as follows, which is a
special case of~\cite[Definition 4.2]{BergmannPrestin2014dlVP} setting~\(l=0\)
therein.

\begin{definition}
  \label{def:dlVP}
  Let \(\mat{M}\in\mathbb Z^{d\times d}\) be a regular matrix
	and~\(g\) be an admissible function.
	The function \(f_{\mat{M},g}\), which is defined by their Fourier
	coefficients as
	\[
		c_{\vect{k}}(f_{\mat{M},g}) \coloneqq \frac{1}{\sqrt{m}}g(\mat{M}^{-\tT}\vect{k}),
		\qquad \vect{k}\in\mathbb Z^d
	\]
	is called \emph{de la Vallée Poussin mean}.
\end{definition}

\begin{figure}[tbp]
	\begin{subfigure}{.5\textwidth}\centering
    \includegraphics[width=.98\textwidth]{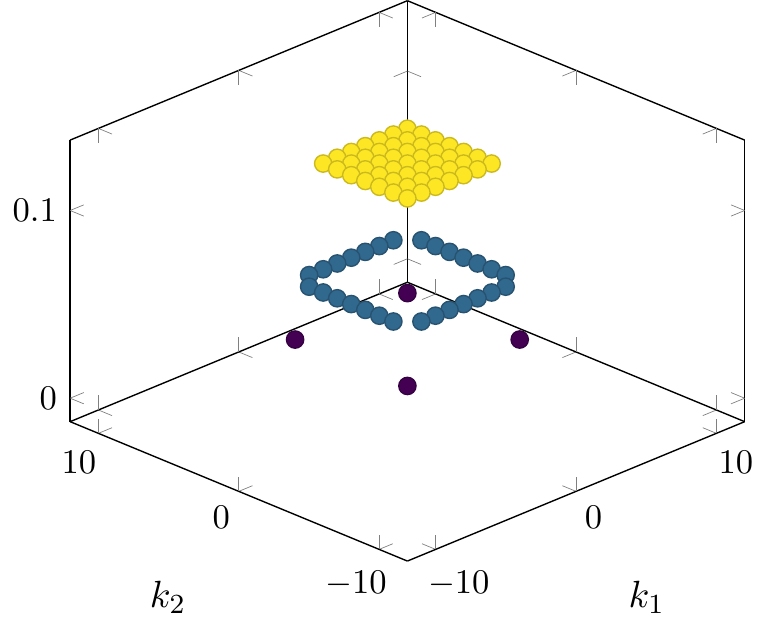}
		\caption{Fourier coefficients of \(f_{\mat{M}_1,\alpha}\),\\
		\(\mat{M}= \bigl(\begin{smallmatrix}
			8&0\\0&8
		\end{smallmatrix}\bigr)\), \(\alpha = \frac{1}{10}\begin{pmatrix}
	1&1
		\end{pmatrix}^\tT\).}
		\label{subfig:dlVP1ck}
	\end{subfigure}
	\begin{subfigure}{.5\textwidth}\centering
    \includegraphics[width=.98\textwidth]{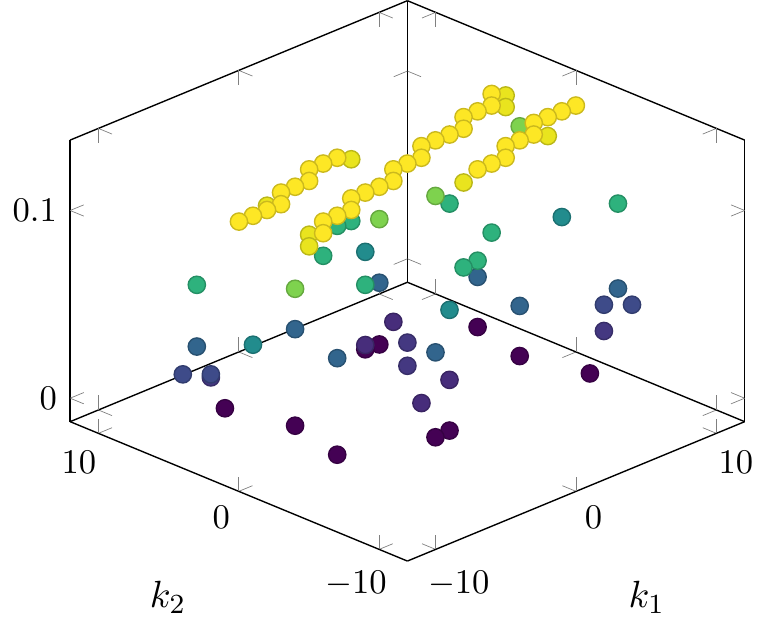}
		\caption{
		Fourier coefficients of \(f_{\mat{M}_2,\alpha}\), \\
				\(\mat{M}_2= \bigl(\begin{smallmatrix}
			4&-2\\4&14
		\end{smallmatrix}\bigr)\),
    \(\alpha = \frac{1}{10}\begin{pmatrix}
			1&1
		\end{pmatrix}^\tT\).}
		\label{subfig:dlVP2ck}
	\end{subfigure}
	\begin{subfigure}{.5\textwidth}\centering
    \includegraphics[width=.98\textwidth]{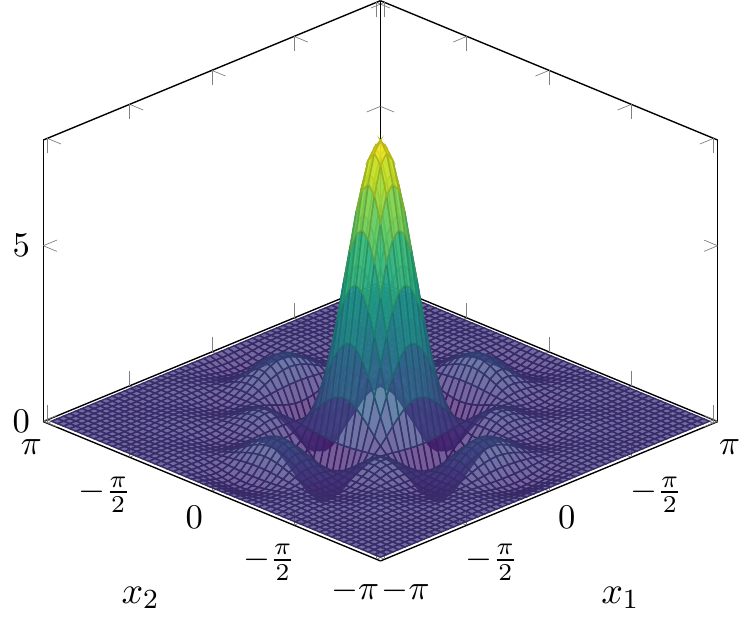}
		\caption{the function \(f_{\mat{M}_1,\alpha}\),\\
		\(\mat{M}= \bigl(\begin{smallmatrix}
			8&0\\0&8
		\end{smallmatrix}\bigr)\), \(\alpha = \frac{1}{10}\begin{pmatrix}
	1&1
		\end{pmatrix}^\tT\).}
		\label{subfig:dlVP1}
	\end{subfigure}
	\begin{subfigure}{.5\textwidth}\centering
    \includegraphics[width=.98\textwidth]{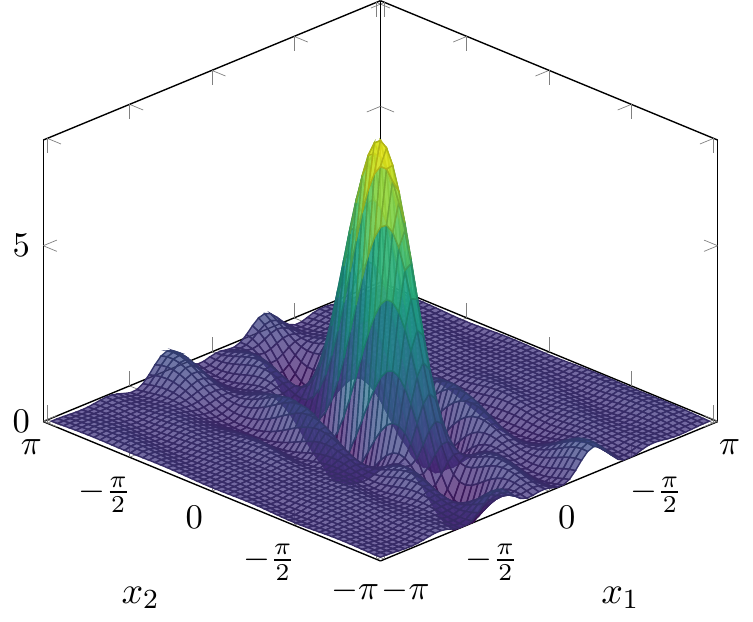}
		\caption{
		The function \(f_{\mat{M}_2,\alpha}\),\\
				\(\mat{M}_2= \bigl(\begin{smallmatrix}
			4&-2\\4&14
		\end{smallmatrix}\bigr)\), \(\alpha = \frac{1}{10}\begin{pmatrix}
			1&1
		\end{pmatrix}^\tT\).}
		\label{subfig:dlVP2}
	\end{subfigure}
	\caption{Two different de la Vallée
	Poussin-means~\(f_{\mat{M}_j,\alpha}\), \(j=1,2\): their Fourier
  coefficients (top row) and their plots in time domain (bottom row).}
	\label{fig:dlVP}
\end{figure}

In case of Box splines \(g_{\vect{\alpha}}\) this generalizes the
one-dimensional de la Vallée Poussin means
from~\cite{Selig:1998,PrestinSelig:1998} to arbitrary patterns including the
tensor product case for diagonal matrices \(\mat{M}\), which where for example used in~\cite{Sprengel1997}. We will use the short hand
notation~\(f_{\mat{M},\vect{\alpha}}\coloneqq f_{\mat{M},g_{\vect{\alpha}}}\)
and omit~\(\alpha\) whenever its clear from the context. It is easy to see,
that by admissibility of~\(g\) the fundamental interpolant exists for any de la
Vallée Poussin mean~\(f_{\mat{M},g}\). The functions \(f_{\mat{M},\alpha}\)
generalize the classical de la Vallée Poussin means to higher dimensions and
anisotropic patterns, for which examples are shown in Figure~\ref{fig:dlVP} and
explained in the following.

Finally, the \Def{Dirichlet kernel} $D_{\mat{M}}$ is defined by using the function $g(\vect{x})$ with
\begin{equation*}
  g(\vect{x}) \coloneqq
  \begin{cases}
    1,   &\vect{x} \in \bigl[ -\frac12, \frac12 \bigr)^d\\
    0,   &\text{otherwise}.
  \end{cases}
\end{equation*}
This kernel is comprised in the definition of the generalized de la Vall\'ee
Poussin mean as well. Furthermore we obtain
the~\emph{modified Dirichlet Kernel} \(f_{\mat{M},\vect{0}}\) as a limiting case
of the de la Vallée Poussin case.

As an example we choose~\(\mat{M}_1  = \bigl( \begin{smallmatrix}8&0\\0&8\end{smallmatrix}\bigr)\) and \(\mat{M}_2= \bigl(\begin{smallmatrix} 4&-2\\4&14\end{smallmatrix}\bigr)\). For \(\mat{M}_1\)
we obtain the usual rectangular (pixel grid) pattern~\(\Pattern(\mat{M}_1)\)
while \(\Pattern(\mat{M}_2)\) models a certain anisotropy, cf.~\cite[Fig.~2.1]{BergmannMerkert2016}. By further setting \(\alpha = \frac{1}{10}\begin{pmatrix}
	1&1
\end{pmatrix}^\tT\) we obtain the de la Vallée Poussin means \(f_{\mat{M}_1,\alpha}\) and \(f_{\mat{M}_2,\alpha}\). Their Fourier
coefficients~\(c_{\vect{k}}(f_{\mat{M}_1,\alpha})\)
and~\(c_{\vect{k}}(f_{\mat{M}_2,\alpha})\) after orthonormalising the
translates, cf.~Lemma~\ref{lem:FI:Props}~\ref{lem:FI:Props:ONB}, are shown in
Figs.~\ref{fig:dlVP}\,\subref{subfig:dlVP1ck} and \subref{subfig:dlVP2ck},
respectively.
Note that the first results in \(64\) translates, while the second determinant
is smaller and results in \(58\) translates.
The functions in time domain are plotted in
Figs.~\ref{fig:dlVP}\,\subref{subfig:dlVP1} and \subref{subfig:dlVP2},
respectively.
While the first can also be obtained by a tensor product of one-dimensional
de la Vallée Poussin means, cf., e.g.~\cite{Selig:1998}, the second one prefers
in time domain certain directions due to its anisotropic form.

\section{Homogenization on spaces of translates}
\label{sec:homogenization}

In the following steps we use anisotropic spaces of translates to discretise the
quasi-static equation of linear elasticity in homogenization. First we introduce
the necessary spaces and differential operators. With these we can state the
partial differential equation we are interested in and two equivalent
formulations, a variational equation and the so-called Lippmann-Schwinger
equation. These formulations make use of the Green operator $\Gamma^0$.
Based on this operator we introduce the periodised Green operator $\Gamma^\p$
and subsequently analyse its properties and special cases. Next we use this
operator to discretise the partial differential equation while splitting the
derivation into two steps.

\subsection{The elasticity problem in periodic homogenization}

FFT-based methods for the equations of linear elasticity in homogenization based
on the Lippmann-Schwinger equation are first introduced by Moulinec and
Suquet~\cite{MoulinecSuquet1994,MoulinecSuquet1998}. Based on their method,
Vond\v{r}ejc~et.~al.~\cite{Vondrejc2014} interpret the resulting discretisation
as a Galerkin projection using trigonometric sums.

In the following we generalize the interpretation using trigonometric sums to
spaces of translates of a periodic function. Therein the trigonometric sums will
appear as a special case, namely when choosing the Dirichlet kernel's
translates and a diagonal matrix $\mat{M}$.

Let $S$ be an arbitrary set, then we introduce the notations
\begin{align*}
  \Sym_d\bigl( S \bigr) &\coloneqq 
  \bigl\lbrace 
  s \in S^{d \times d}: s_{ij} = s_{ji} \in S \text{ for all } i,j=1,\dots,d
  \bigr\rbrace
  \label{eq:Symd},\\
  \SSym_d\bigl( S \bigr) &\coloneqq 
  \bigl\lbrace 
  s \in \Sym_d(S) \times \Sym_d(S): s_{ijkl} = s_{klij} \in S\\
  &\qquad\qquad\qquad\qquad\qquad\qquad
    \text{ for all }i,j,k,l=1,\dots,d\bigr\rbrace.
\end{align*}
The space $\Sym_d(S)$ corresponds to symmetric matrices built of elements of $S$
and $\SSym_d(S)$ corresponds to fourth-order tensors $\Stiffness = \bigl(
\Stiffness_{ijkl} \bigr)_{ijkl} \in S^{d \times d \times d \times d}$ with minor
and major symmetries, i.e. $\Stiffness_{ijkl} = \Stiffness_{jikl} =
\Stiffness_{ijlk} = \Stiffness_{klij}$.

We endow the space of symmetric matrices $\Sym_d\bigl(\R\bigr)$ with the
Frobenius inner product
\begin{equation*}
  \langle A, B \rangle \coloneqq \sum_{i,j=1}^d
  A_{ij}B_{ij},\quad A,B \in \Sym_d\bigl(\R\bigr).
\end{equation*}

We call a function $\mathcal{A} \in \SSym_d\bigl(L^2\bigl(\T^d\bigr)\bigr)$
\Def{uniformly elliptic} if there exist constants $0 < l_{\mathcal{A}^0} \leq
u_{\mathcal{A}^0} < \infty$ such that for almost all functions $\gamma \in
\Sym_d\bigl(L^2\bigl(\T^d\bigr)\bigr)$ it holds true that
\begin{equation*}
  l_{\mathcal{A}^0} \norm{\gamma}^2 \leq
  \bigl\langle
  \mathcal{A} \gamma, \gamma
  \bigr\rangle
  \leq u_{\mathcal{A}^0} \norm{\gamma}^2.
\end{equation*}
A uniformly elliptic and constant function $\mathcal{A}^0 \in
\SSym_d\bigl(L^2\bigl(\T^d\bigr)\bigr)$ can be identified with an element from
$\SSym_d\bigl(\R\bigr)$ and is called \Def{elliptic} in the following.

The Sobolev space $H^1(\T^d)$ is defined via Fourier series, i.e.
\begin{align*}
  H^1\bigl(\T^d\bigr) &\coloneqq 
  \bigl\lbrace
  f \in L^1(\T^d): \norm{f}[H^1] < \infty
  \bigr\rbrace
  \intertext{with the norm}
  \norm{f}[H^1] &\coloneqq 
  \biggl\lVert
  \sum_{\vect{k} \in \Z^d} 
  \bigl(
    1 + \norm{k}[2]^2
  \bigr)^{\frac12} 
  c_\vect{k}(f) e^{\im \vect{k}^\tT \cdot}
  \biggr\rVert_{L^2}.
\end{align*}

To simplify the equations of linear elasticity and the differential operators
occurring therein we additionally introduce for $\vect{u} \in H^1(\T^d)^d$ the
\Def{symmetrized gradient operator}
\begin{equation}
  \GradSym(\vect{u}) \coloneqq \frac12 \Bigl( \nabla \vect{u} + \bigl( \nabla
  \vect{u} \bigr)^\tT \Bigr) \in \Sym_d\bigl(L^2\bigl(\T^d\bigr)\bigr),
  \label{eq:GradSym}
\end{equation}
and the \Def{divergence operator} $\Div(\vect{u})$ as the formal $L^2$-adjoint of
$\GradSym(\vect{u})$. The action of the symmetrized gradient operator in Fourier
space is for $\vect{k} \in \Z^d$ given by
\begin{equation*}
    c_\vect{k}\bigl(\GradSym(\vect{u})\bigr) =
    \GradSym_\vect{k} c_\vect{k}(\vect{u}) \coloneqq
      \frac{\im}{2} \bigl(\vect{k} c_\vect{k}(\vect{u})^\tT +
      c_\vect{k}(\vect{u}) \vect{k}^\tT\bigr).
\end{equation*}

The basic solution space for the PDE we want to analyse is given by \Def{symmetric
gradient fields with zero mean}
\begin{equation*}
  \Curlfree\bigl(\T^d\bigr) \coloneqq
  \Bigl\lbrace
  \varepsilon \in \Sym_d\bigl(L^2\bigl(\T^d\bigr)\bigr):
  \exists \vect{u} \in H^1\bigl(\T^d\bigr)^d, \varepsilon = \GradSym \vect{u}
  \Bigr\rbrace.
\end{equation*}

With these preparations we can now state the partial differential equation
describing \Def{quasi-static linear elasticity in homogenization}.
\begin{definition}
  \label{def:weak_pde}
  Let $\Stiffness \in
  \SSym_d\bigl(L^{\infty}\bigl(\T^d\bigr)\bigr)$ be uniformly
  elliptic and let $\varepsilon^0 \in
  \Sym_d\bigl(\R\bigr)$ be given. We want to find the strain
  $\varepsilon \in \Curlfree\bigl(\T^d\bigr)$ such that
  \begin{equation}\label{eq:weak_pde}
    \bigl\langle
    \Stiffness \tensorProd \varepsilon,\tilde{\gamma}
    \bigr\rangle = 
    -\bigl\langle
    \Stiffness
    \tensorProd \varepsilon^0,\tilde{\gamma}
    \bigr\rangle
  \end{equation}
  holds for all $\tilde{\gamma} \in \Curlfree\bigl(\T^d\bigr)$. With $\Stiffness \tensorProd
  \varepsilon$ we hereby denote the product of the fourth-order stiffness tensor
  $\Stiffness$ and the second-order strain $\varepsilon$, the symmetry making the
  order of multiplication irrelevant.
\end{definition}
The stiffness distribution $\Stiffness$ describes the material behaviour and is
in practical applications usually a piece-wise constant function. The role of
the macroscopic strain $\varepsilon^0$ is that of an overall strain that is applied
to the composed material and corresponds to pulling (or compressing) the
composite in a certain direction.

For applications one is interested either in the strain field $\varepsilon$
---or derivates like stress or displacements, respectively--- or in the
so-called \Def{effective stiffness} of the medium. The overall
stiffness behaviour of the composite given by the action of $\Stiffness^\eff \in
\SSym_d\bigl(\R\bigr)$ on $\varepsilon^0 \in \Sym_d\bigl(\R\bigr)$ is defined by
\begin{equation*}
  \Stiffness^\eff \tensorProd \varepsilon^0 \coloneqq \frac{1}{(2 \pi)^d}
  \int_{\T^d} \Stiffness(x) \tensorProd \varepsilon(x) \dx x
\end{equation*}
where $\varepsilon$ is the solution of~\eqref{eq:weak_pde} corresponding to
$\varepsilon^0$.

The space of suitable ansatz functions $\Curlfree\bigl(\T^d\bigr)$ is in
practice difficult to work with, especially regarding the discretisation steps
that follow. A method to deal with that is introduced by
Vond\v{r}ejc~et.~al.~\cite{Vondrejc2014}. They derive a projection operator
$\Gamma^0 \Stiffness^0$ with $L^2$-adjoint $\Stiffness^0 \Gamma^0$ that maps
$\Sym_d\bigl(L^2\bigl(\T^d\bigr)\bigr)$ onto $\Curlfree\bigl(\T^d\bigr)$ and
thus replaces the structurally complicated space by a simpler one necessitating
a more involved PDE. The operator $\Gamma^0$ acts as a second order derivative
of a preconditioner that solves the constant coefficient PDE $\Div \Stiffness^0
\GradSym \vect{u} = \vect{f}$. For a derivation of the formula see for
example~\cite{MoulinecSuquet1994,MoulinecSuquet1998,Schneider2016}.

\begin{definition}
  \label{def:projection}
  For a constant elliptic reference stiffness $\Stiffness^0 \in
  \SSym_d\bigl(\R\bigr)$ define the \Def{Green operator} $\Gamma^0 \colon
  \Sym_d\bigl(L^2\bigl(\T^d\bigr)\bigr) \rightarrow \Curlfree(\T^d)$ which acts
  as a Fourier multiplier. Its action on a field $\varepsilon \in
  \Sym_d\bigl(L^2\bigl(\T^d\bigr)\bigr)$ is given by the equation
  \begin{align}\label{eq:gamma}
      \Gamma^0 \tensorProd \varepsilon &\coloneqq
      \sum_{\vect{k} \in \Z^d}
      \hat{\Gamma}^0_\vect{k} \tensorProd c_\vect{k}(\varepsilon)
      e^{2 \pi \im \vect{k}^\tT \cdot}
      \intertext{with equality in $L^2$-sense and Fourier coefficients}
      \nonumber
      \hat{\Gamma}^0_\vect{k} \tensorProd c_\vect{k}(\varepsilon)
      &\coloneqq
      \GradSym_\vect{k}
      \Bigl(
	\overline{\GradSym_\vect{k}}^\tT \tensorProd \Stiffness^0 \tensorProd
	\GradSym_\vect{k}
      \Bigr)^{-1}
      \overline{\GradSym_\vect{k}}^\tT c_\vect{k}(\varepsilon),\quad
      \vect{k} \in\Z^d.
  \end{align}
\end{definition}

This operator allows to reformulate~\eqref{eq:weak_pde} with test functions in
$\Sym_d\bigl(L^2\bigl(\T^d\bigr)\bigr)$ by projecting them onto the required
function space. Further, this additional operator is easy do apply as it acts as
a convolution operator, i.e. a Fourier multiplier.
\begin{proposition}
  \label{def:weak_pde_projected}
  Let $\varepsilon^0 \in \Sym_d\bigl(\R\bigr)$, then $\varepsilon \in
  \Curlfree\bigl(\T^d\bigr)$ fulfils
  \begin{equation*}
    \bigl\langle
    \Stiffness \tensorProd \varepsilon,\tilde\gamma
    \bigr\rangle =
    -\bigl\langle
    \Stiffness \tensorProd \varepsilon^0,\tilde\gamma
    \bigr\rangle
  \end{equation*}
  for all $\tilde{\gamma} \in \Curlfree\bigl(\T^d\bigr)$ if and only if
  \begin{equation}\label{eq:weak_pde_projected}
    \bigl\langle
    \Stiffness^0 \Gamma^0 \Stiffness \tensorProd \varepsilon,\gamma
    \bigr\rangle =
    -\bigl\langle
    \Stiffness^0 \Gamma^0 \Stiffness \tensorProd \varepsilon^0,\gamma
    \bigr\rangle
  \end{equation}
  holds true for all $\gamma \in \Sym_d\bigl(L^2\bigl(\T^d\bigr)\bigr)$.
  \begin{proof}
    For a proof see~\cite[Proposition 3]{Vondrejc2014}.
  \end{proof}
\end{proposition}

By the properties of the operator $\Gamma^0 \Stiffness^0$ this equation is equivalent to the so-called
\Def{Lippmann-Schwinger equation} with fixed point $\varepsilon$, cf.\cite[Lemma 2]{Vondrejc2014}:
\begin{proposition}
  \label{cor:weak_pde_projected_2}
  Let $\varepsilon^0 \in \Sym_d\bigl(\R\bigr)$, then $\varepsilon \in
  \Curlfree\bigl(\T^d\bigr)$ fulfils
  \begin{equation*}
    \bigl\langle
    \Stiffness \tensorProd \varepsilon,\tilde\gamma
    \bigr\rangle =
    -\bigl\langle
    \Stiffness \tensorProd \varepsilon^0,\tilde\gamma
    \bigr\rangle
  \end{equation*}
  for all $\tilde{\gamma} \in \Curlfree\bigl(\T^d\bigr)$ if and only if
  \begin{equation}\label{eq:lippmann-schwinger}
    \bigl\langle
    \varepsilon +
    \Gamma^0 \bigl(
    \Stiffness - \Stiffness^0
  \bigr) \tensorProd
  \bigl(
  \varepsilon + \varepsilon^0
\bigr),
\gamma
\bigr\rangle = 0
  \end{equation}
is fulfilled for all $\gamma \in \Sym_d\bigl(L^2\bigl(\T^d\bigr)\bigr)$.
\end{proposition}

\begin{remark}
  \label{rem:Gamma_identity}
  The proof makes use of the adjoint operator $\Stiffness^0 \Gamma^0$ and the
  identity
\begin{equation}
  \varepsilon - \Gamma^0 \Stiffness^0 \tensorProd \varepsilon = \varepsilon^0
  \label{eq:projection_identity}
\end{equation}
which holds in weak sense. The above identity is equivalent to $\Gamma^0
\Stiffness^0$ being a projection operator that maps constants to zero. These
properties are proven in~\cite[Lemma 2]{Vondrejc2014}.
\end{remark}

With this at hand we can now proceed to discretising the Lippmann-Schwinger
equation in a space of translates~\(V_{\mat{M}}^f\), \(f\in A(\mathbb T^d)\).

\subsection{The periodised Green operator}

The Lippmann-Schwinger equation is first discretised by Moulinec and
Suquet~\cite{MoulinecSuquet1994,MoulinecSuquet1998} using a Fourier collocation
scheme, cf.~\cite{ZVNM:2010NumHom}. The resulting fixed point algorithm
inspired many publications, for an overview see for example~\cite{Mishra2016}.

In contrast, Vond\v{r}ejc~et.~al.~\cite{Vondrejc2014} employ a Galerkin projection
of~\eqref{eq:weak_pde_projected} using trigonometric polynomials as ansatz
functions and obtain the same discretisation. In the following we want to
generalize this approach to spaces of translates on anisotropic lattices.

Throughout this section $\mat{M} \in \Z^{d \times d}$ denotes a regular pattern
matrix that defines a translation invariant space $V_\mat{M}^f$ spanned by the
translates $\Translate(\vect{y}) f$, $\vect{y} \in \Pattern(\mat{M})$, such
that a fundamental interpolant exists,
see Lemma~\ref{lem:FI:Props}~\ref{lem:FI:Prosp:Existence}. Especially, this implies that the translates
$\Translate(\vect{y}) f$ of $f$ are linearly independent.

This fundamental interpolant is denoted by
$\Fundamental_\mat{M} \in V_\mat{M}^f$, i.e. there exist coefficients
$\hat{a}_\vect{h}$, $\vect{h} \in \generatingSet(\mat{M}^\tT)$ such that
\begin{equation}
  \label{eq:fundamental_coeff}
  c_{\vect{h}+\mat{M}^\tT \vect{z}}(\Fundamental_{\mat{M}}) =
  \hat{a}_\vect{h} c_{\vect{h}+\mat{M}^\tT \vect{z}}(f)
\end{equation}
holds true for all $\vect{h} \in \generatingSet(\mat{M}^\tT)$ and
$\vect{z} \in \Z^d$.

From now on we assume for functions $\gamma
\in \Sym_d\bigl(V^f_\mat{M} \bigr)$ that
\begin{equation*}
  \gamma = \sum_{\vect{y} \in \Pattern(\mat{M})} 
  G_\vect{y} \Translate(\vect{y}) f
\end{equation*}
holds true. We further denote the discrete Fourier transform
of~\(\vect{G}
= (G_{\vect{y}})_{\vect{y}\in\Pattern(\mat{M})}\)
by~\(\vect{\hat G} = (\hat G_{\vect{h}})_{\vect{h}\in\generatingSet(\mat{M}^\tT)} = \Fourier(\mat{M})\vect{G}\)

A Galerkin projection of~\eqref{eq:weak_pde_projected} onto the space of
translates requires the definition of a Green operator similar to
Definition~\ref{def:projection}. To account for the finite dimensional space
$V_\mat{M}^f$ the operator $\Gamma^0$ has to be periodised in frequency domain.

\begin{definition}
  \label{def:gamma_tilde_fourier}
  We call the Fourier multiplier $\Gamma^\p$ the \Def{periodised Green
  operator on $V^f_\mat{M}$} and define its action onto a field $\gamma \in
  \Sym_d\bigl(V^f_\mat{M}\bigr)$
  by
  \begin{align}
    \label{eq:gamma_tilde_fourier}
    \Gamma^\p \tensorProd \gamma
    &\coloneqq
    \sum_{\vect{y} \in \Pattern(\mat{M})} 
    \Gamma^\p_\vect{y} \tensorProd G_\vect{y} \Translate(\vect{y}) f.
    \intertext{In terms of Fourier sums this is the same as}
    \nonumber
    \Gamma^\p \tensorProd \gamma
    &\coloneqq
    \sum_{\vect{h} \in \generatingSet(\mat{M}^\tT)}
    \hat{\Gamma}^\p_\vect{h} \tensorProd \hat{G}_\vect{h}
    c_\vect{h}^\mat{M}(f) e^{2 \pi \im \vect{h}^\tT \cdot}
    \intertext{with Fourier coefficients}
    \nonumber
    \hat{\Gamma}^\p_\vect{h} \tensorProd \hat{G}_\vect{h}
    &\coloneqq
    m \Bigl[
      \bigl\lbrace
      \hat{\Gamma}^0_\vect{k} \abs{c_\vect{k}(f)}^2
      \bigr\rbrace_{\vect{k} \in \Z^d}
    \Bigr]^\mat{M}_\vect{h} \tensorProd
    \hat{G}_\vect{h},\ \vect{h} \in \generatingSet(\mat{M}^\tT).
  \end{align}
\end{definition}
%
%
\paragraph{Properties of the periodised Green operator.}
In the trigonometric collocation case of Vond\v{r}ejc~et.al.~\cite{Vondrejc2014}
the Green operator keeps the same form after discretisation, i.e.~a restriction
of the Fourier series to a bounded cube. This also holds true for the
generalization to anisotropic patterns~\cite{BergmannMerkert2016} where the cube
is replaced by a parallelotope, i.e.~to the set $\generatingSet(\mat{M}^\tT)$.
The properties of the Green operator and the projection operator $\Gamma^0
\Stiffness^0$ are shown via properties of its Fourier coefficients. Hence, the
proofs in the continuous and the discretised case can be done analogously. This
is no longer the case for the approach using translation invariant spaces.

\begin{theorem}
  The operator $\Gamma^\p \Stiffness^0$ has the $L^2$-adjoint $\Stiffness^0
  \Gamma^p$.
  \begin{proof}
    The proof for $\Gamma^0 \Stiffness^0$ in~\cite[Lemma 2 (ii)]{Vondrejc2014}
    relies purely on the symmetry of the operator. This symmetry is preserved by
    the periodised Green operator~\(\Gamma^\p \Stiffness^0\) and therefore the proof is analogous.
  \end{proof}
\end{theorem}

Vond\v{r}ejc et.~al.~introduce the operator $\Gamma^0 \Stiffness^0$ to project
functions in $\Sym_d\bigl(L^2\bigl(\T^d\bigr)\bigr)$ onto
$\Curlfree\bigl(\T^d\bigr)$. The operator $\Gamma^\p \Stiffness^0$ has similar
properties and maps onto the respective discretised versions of these spaces.

\begin{theorem}
  For all $\gamma \in \Sym_d \bigl(V_\mat{M}^f \bigr)$ is holds true that $\Gamma^\p
  \Stiffness^0 \tensorProd \gamma \in \Curlfree\bigl(\T^d\bigr) \cap
  \Sym_d\bigl(V_\mat{M}^f\bigr)$.
  \begin{proof}
    For $\vect{y} \in \Pattern(\mat{M})$ we have that
    \begin{align*}
      \bigl(\Gamma^\p \Stiffness^0 \tensorProd \gamma\bigr)(\vect{y}) &=
      \sum_{\vect{h} \in \generatingSet(\mat{M}^\tT)} 
      \hat{\Gamma}^\p_\vect{h}\Stiffness^0  \tensorProd
      \hat G_\vect{h} c_\vect{h}^\mat{M}(f)
      e^{2 \pi \im \vect{h}^\tT \vect{y}}\\
      &= \sum_{\vect{h} \in \generatingSet(\mat{M}^\tT)} 
      m\bigl[
	\bigl\lbrace
	  \hat{\Gamma}^0_\vect{k} \abs{c_\vect{k}(f)}^2
	  \bigr\rbrace_{\vect{k} \in \Z^d}
	\bigr]_\vect{h}^\mat{M}
	\Stiffness^0 \tensorProd
      \hat G_\vect{h} c_\vect{h}^\mat{M}(f)
      e^{2 \pi \im \vect{h}^\tT \vect{y}}.
    \end{align*}
    First,  observe this can be rewritten with
    \begin{equation*}
      \tilde{G}_\vect{h} \coloneqq
      m\bigl[
	\bigl\lbrace
	  \hat{\Gamma}^0_\vect{k} \abs{c_\vect{k}(f)}^2
	  \bigr\rbrace_{\vect{k} \in \Z^d}
	\bigr]_\vect{h}^\mat{M}
	\Stiffness^0 \tensorProd
      \hat G_\vect{h}
    \end{equation*}
    for $\vect{h} \in \generatingSet(\mat{M}^\tT)$ as
    \begin{equation*}
      \bigl(\Gamma^\p \Stiffness^0 \tensorProd \gamma\bigr)(\vect{y})
      = \sum_{\vect{h} \in \generatingSet(\mat{M}^\tT)} 
      \tilde{G}_\vect{h} c_\vect{h}^\mat{M}(f)
      e^{2 \pi \im \vect{h}^\tT \vect{y}}.
    \end{equation*}
    Using~\eqref{eq:inTranslates:ck} the result is a function in $\Sym_d\bigl( V_\mat{M}^f \bigr)$.

    Expanding the bracket sum and the Green operator $\Gamma^0$ and yields
    \begin{align*}
      \bigl(\Gamma^\p \Stiffness^0 \tensorProd \gamma\bigr)(\vect{y})
      &= \sum_{\vect{h} \in \generatingSet(\mat{M}^\tT)} 
      m \sum_{\vect{z} \in \Z^d}
      \abs{c_{\vect{h}+\mat{M}^\tT \vect{z}}(f)}^2
      \hat{\Gamma}^0_{\vect{h}+\mat{M}^\tT \vect{z}}
      \Stiffness^0 \tensorProd
      \hat G_\vect{h} c_\vect{h}^\mat{M}(f)
      e^{2 \pi \im \vect{h}^\tT \vect{y}}\\
      &= \sum_{\vect{h} \in \generatingSet(\mat{M}^\tT)} 
      m \sum_{\vect{z} \in \Z^d}
      \abs{c_{\vect{h}+\mat{M}^\tT \vect{z}}(f)}^2
      \GradSym_{\vect{h}+\mat{M}^\tT \vect{z}}\\
      &\qquad \times \bigl( 
	\overline{\GradSym_{\vect{h}+\mat{M}^\tT \vect{z}}}^\tT \tensorProd
	\Stiffness^0 \tensorProd
	\GradSym_{\vect{h}+\mat{M}^\tT \vect{z}}
      \bigr)^{-1}
      \overline{\GradSym_{\vect{h}+\mat{M}^\tT \vect{z}}}^\tT \\
      &\qquad \times \Stiffness^0 \tensorProd \hat G_\vect{h} c_\vect{h}^\mat{M}(f)
      e^{2 \pi \im \vect{h}^\tT \vect{y}}.
    \end{align*}
    We define new Fourier coefficients
    \begin{align*}
      c_{\vect{h}+\mat{M}^\tT \vect{z}}(\vect{\tilde{u}})
      \coloneqq 
      &\,m 
      \abs{c_{\vect{h}+\mat{M}^\tT \vect{z}}(f)}^2
      \bigl( 
	\overline{\GradSym_{\vect{h}+\mat{M}^\tT \vect{z}}}^\tT \tensorProd
	\Stiffness^0 \tensorProd
	\GradSym_{\vect{h}+\mat{M}^\tT \vect{z}}
      \bigr)^{-1}
      \overline{\GradSym_{\vect{h}+\mat{M}^\tT \vect{z}}}^\tT \\
      &\quad \times\Stiffness^0 \tensorProd
      \hat G_\vect{h} c_\vect{h}^\mat{M}(f).
    \end{align*}
    With the decomposition $\vect{k} = \vect{h} + \mat{M}^\tT \vect{z}$
    for $\vect{h}\in \generatingSet(\mat{M}^\tT)$ and $\vect{z} \in \Z^d$ the Fourier
    coefficients from the formula above can be collected with respect to congruence classes of the generating set~\(\generatingSet(\mat{M}^\tT\)
    and this yields
    \begin{equation*}
      \bigl(\Gamma^\p \Stiffness^0 \tensorProd \gamma\bigr)(\vect{y})
      = \sum_{\vect{k} \in \Z^d}
      \GradSym_{k} c_\vect{k}(\vect{\tilde{u}}) 
      e^{2 \pi \im \vect{k}^\tT \vect{y}}
    \end{equation*}
    for $\vect{y} \in \Pattern(\mat{M})$.

    We finally take a closer look at the Fourier series
    \begin{equation}\label{eq:GradSym_series}
     \sum_{\vect{k} \in \Z^d}
     c_\vect{k}(\vect{\tilde{u}})
      e^{2 \pi \im \vect{k}^\tT \vect{y}}
    \end{equation}
    and analyse its convergence. With $f \in A(\T^d) \subset L^2(\T^d)$ and
    because the Fourier coefficients $c_\vect{h}(\gamma)$ and
    $c_\vect{h}^\mat{M}(f)$ depend only on $\vect{h} \in
    \generatingSet(\mat{M}^\tT)$ and not on $\vect{z} \in \Z^d$ the Fourier series
    \begin{equation*}
    \sum_{\vect{h} \in \generatingSet(\mat{M}^\tT)} \sum_{\vect{z} \in \Z^d}
      m \abs{ c_{\vect{h} + \mat{M}^\tT \vect{z}}(f) }^2
      \Stiffness^0 \tensorProd
      \hat G_\vect{h} c_\vect{h}^\mat{M}(f)
      e^{2 \pi \im (\vect{h}+\mat{M}^\tT \vect{z})^\tT \vect{y}}
    \end{equation*}
    converges and the result is at least in
    $\Sym_d\bigl(L^2\bigl(\T^d\bigr)\bigr)$.
    When we apply the differential operators $\overline{\GradSym_\vect{k}}^\tT$ and
    $\bigl(\overline{\GradSym_\vect{k}}^\tT \Stiffness^0 \GradSym_\vect{k}
    \bigr)^{-1}$ we differentiate the $L^2$ function once and integrate twice.
    The resulting function of interest is at least once weakly differentiable,
    i.e. in $H^1(\T^d)^d$. Thus, another application of
    $\GradSym$ is admissible and the Fourier series~\eqref{eq:GradSym_series}
    converges. Further~\eqref{eq:GradSym_series} is a gradient field with mean
    zero, cf.~\eqref{eq:GradSym}, and the proof is concluded.
  \end{proof}
\end{theorem}

A special choice for the space $V^f_\mat{M}$ comes from using the Dirichlet
kernel~\(D_{\mat{M}}\), where $\Gamma^\p$ coincides with $\Gamma^0$, which also
occurs in the derivation in~\cite{MoulinecSuquet1998}.

\begin{theorem}\label{thm:CoincidenceGamma0GammaP}
  For the Dirichlet kernel $D_{\mat{M}}$ the periodised Green operator $\Gamma^\p$
  of~\eqref{eq:gamma_tilde_fourier} on $V^{D_{\mat{M}}}_\mat{M}$ coincides with
  the Green operator $\Gamma^0$.
  \begin{proof}
    When we insert the formula for the Dirichlet kernel $D_{\mat{M}}$
    into~\eqref{eq:gamma_tilde_fourier}, the sum reduces to one single term
    which is exactly $\hat{\Gamma}^0_\vect{k}$ and thus the proof is completed.
  \end{proof}
\end{theorem}

In contrast to the operator $\Gamma^0 \Stiffness^0$, the periodised Green
operator $\Gamma^\p \Stiffness^0$ corresponding to~\(V_{\mat{M}}^f\) is in general no longer a projection. 
\begin{theorem}
  The periodised Green operator $\Gamma^\p \Stiffness^0$ corresponding to \(V_{\mat{M}}^f\) is a projection
  operator if and only if $V_{\mat{M}}^f=V_{\mat{M}}^{D_{\mat{M}}}$, i.e. iff either \(f\) or (one of) its orthonormalised translates is the Dirichlet kernel~\(D_{\mat{M}}\).
  \begin{proof}
    Consider for a field $\gamma \in \Sym_d\bigl(V_\mat{M}^f\bigr)$ the Fourier
    series
    \begin{equation*}
      \sum_{\vect{h} \in \generatingSet(\mat{M}^\tT)} 
      \hat{\Gamma}^\p_\vect{h} \Stiffness^0
      \hat{\Gamma}^\p_\vect{h} \Stiffness^0
      \hat G_\vect{h} c_\vect{h}^\mat{M}(f) 
      e^{2 \pi \im \vect{h}^\tT \vect{y}}
    \end{equation*}
    for $\vect{y} \in \Pattern(\mat{M})$ and insert the definition of the
    periodised Green operator~\eqref{eq:gamma_tilde_fourier} to get
    \begin{equation}\label{eq:Gamma_projection}
      \begin{split}
      \sum_{\vect{h} \in \generatingSet(\mat{M}^\tT)} 
      \sum_{\vect{z},\vect{z}' \in \Z^d}
      m^2&
      \abs{c_{\vect{h} + \mat{M}^\tT \vect{z}}(f)}^2
      \abs{c_{\vect{h} + \mat{M}^\tT \vect{z}'}(f)}^2
      \\
      &\times\ 
      \hat{\Gamma}^0_{\vect{h} + \mat{M}^\tT \vect{z}} 
      \Stiffness^0 
      \hat{\Gamma}^0_{\vect{h} + \mat{M}^\tT \vect{z}'} 
      \Stiffness^0
      \hat G_\vect{h} c_\vect{h}^\mat{M}(f) 
      e^{2 \pi \im \vect{h}^\tT \vect{y}}.
      \end{split}
    \end{equation}
    From~\cite[Lemma 2 (iii)]{Vondrejc2014} we know that for $\vect{k} \in \Z^d$
    we have that $\hat{\Gamma}^0_\vect{k} \Stiffness^0 \hat{\Gamma}^0_\vect{k}
    \Stiffness^0 = \hat{\Gamma}^0_\vect{k} \Stiffness^0$, i.e. that
    $\Stiffness^0 \Gamma^0$ is a projection. This does not hold true for the
    mixed terms in~\eqref{eq:Gamma_projection}, i.e. summands with $\vect{z} \neq
    \vect{z}'$. These only vanish if $c_{\vect{k}}(f) = 0$ for $\vect{k} \not\in
    \generatingSet(\mat{M}^\tT)$, i.e. \(f\in V_{\mat{M}}^{D_{\mat{M}}}\),
    cf.~Theorem~\ref{thm:CoincidenceGamma0GammaP}
  \end{proof}
\end{theorem}
Hence \(\Gamma^\p\) is a projection if and only if $\Gamma^\p = \Gamma^0$ with
respect to the corresponding generating set~\(\generatingSet(\mat{M}^\tT)\). In
addition the periodised Green operator~\(\Gamma^\p\Stiffness^0\) is bounded with
the same bound as~\(\Gamma^0\Stiffness^0\).

\begin{theorem}
  Let the translates of $f$ be orthonormal and let $\Stiffness^0 \in
  \SSym_d\bigl(\R\bigr)$ be elliptic with constants $0 < l_{\Stiffness^0} \leq
  u_{\Stiffness^0} < \infty$. Then the periodised Green operator $\Gamma^\p \Stiffness^0$ corresponding to $V^f_\mat{M}$, is bounded by
  \begin{equation*}
    \norm[\bigl]{\Gamma^p \Stiffness^0 \tensorProd
    \gamma}[\Sym_d(L^2(\T^d))]^2 \leq
	  \frac{u_{\Stiffness^0}}{l_{\Stiffness^0}}
	  \norm{\gamma}[\Sym_d(L^2(\T^d))]
  \end{equation*}
  for all $\gamma \in \Sym_d\bigl( V^f_\mat{M} \bigr)$.

  \begin{proof}
    The Parseval equation together with the splitting $\vect{k} = \vect{h} +
    \mat{M}^\tT \vect{z}$ with $\vect{h} \in \generatingSet(\mat{M}^\tT)$ and
    $\vect{z} \in \Z^d$ yields 
	\begin{align*}
	  \norm[\bigl]{
	    \Gamma^p \Stiffness^0 &\tensorProd \gamma
	    }[\Sym_d(L^2(\T^d))]^2
	  \\
    &=
	  \sum_{\vect{h} \in \generatingSet(\mat{M}^\tT)}
	  \sum_{\vect{z} \in \Z^d}
	  \norm[\bigl]{
	    \hat{\Gamma}^\p_\vect{h} \Stiffness^0 \tensorProd
	    \hat G_\vect{h} c_{\vect{h} + \mat{M}^\tT \vect{z}}(f)
	  }^2\\
	  &=
	  \sum_{\vect{h} \in \generatingSet(\mat{M}^\tT)}
	  \sum_{\vect{z} \in \Z^d}
	  \norm[\biggl]{
	    m
	    \Bigl[
	      \bigl\lbrace
	      \hat{\Gamma}^0_\vect{k} \abs{c_\vect{k}(f)}^2
	      \bigr\rbrace_{\vect{k} \in \Z^d}
	    \Bigr]^\mat{M}_\vect{h}
	     \Stiffness^0 \tensorProd 
	    \hat G_\vect{h} c_{\vect{h} + \mat{M}^\tT \vect{z}}(f)
	   }^2.
	\end{align*}
	The Cauchy-Schwarz theorem together with inserting the formula for the
	bracket sums~\eqref{eq:bracketsum} bounds this expression from above by
	\begin{equation*}
    \begin{split}      
	  \norm[\bigl]{\Gamma^p \Stiffness^0 &\tensorProd \gamma}[\Sym_d(L^2(\T^d))]^2
       \\&\leq
	  \sum_{\vect{h} \in \generatingSet(\mat{M}^\tT)}
	  \sum_{\vect{z} \in \Z^d}
	  m^2
	  \sum_{\vect{z'} \in \Z^d}
	  \norm[\bigl]{\hat{\Gamma}^0_{\vect{h} + \mat{M}^\tT \vect{z'}} \Stiffness^0
	    \tensorProd 
	    \hat G_\vect{h} c_{\vect{h} + \mat{M}^\tT \vect{z}}(f)
	  }
	\abs{c_{\vect{h}+\mat{M}^\tT \vect{z'}}(f)}^4
    \end{split}
	\end{equation*}
	A standard estimate for $\hat{\Gamma}^0_\vect{k}$ with $\vect{k} \in
	\Z^d$ is $\norm{\hat{\Gamma}^0_\vect{k}}^2 \leq
	\frac{1}{l_{\Stiffness^0}}$, see e.g.~\cite{Vondrejc2014}.
	This leads to
	\begin{equation*}
    \begin{split}
  	  \norm[\bigl]{\Gamma^p \Stiffness^0 &\tensorProd \gamma}[\Sym_d(L^2(\T^d))]^2
	  \\&\leq
	  \sum_{\vect{h} \in \generatingSet(\mat{M}^\tT)}
	  \sum_{\vect{z} \in \Z^d}
	  m^2
	  \frac{u_{\Stiffness^0}}{l_{\Stiffness^0}}
	  \norm{\hat G_\vect{h} c_{\vect{h} + \mat{M}^\tT \vect{z}}(f)}^2
	  \sum_{\vect{z'} \in \Z^d}
	  \abs{c_{\vect{h}+\mat{M}^\tT \vect{z'}}(f)}^4
    \end{split}
	\end{equation*}
	and Jensen's inequality together with
	Lemma~\ref{lem:FI:Props}~\ref{lem:FI:Props:ONB}  results in
	\begin{align*}
	  \norm[\bigl]{\Gamma^p \Stiffness^0 &\tensorProd \gamma}[\Sym_d(L^2(\T^d))]^2
	  \\&\leq
	  \sum_{\vect{h} \in \generatingSet(\mat{M}^\tT)}
	  \sum_{\vect{z} \in \Z^d}
	  m^2
	  \frac{u_{\Stiffness^0}}{l_{\Stiffness^0}}
	  \norm{\hat G_\vect{h} c_{\vect{h} + \mat{M}^\tT \vect{z}}(f)}^2
	  \bigl(
	    \sum_{\vect{z'} \in \Z^d}
	  \abs{c_{\vect{h}+\mat{M}^\tT \vect{z'}}(f)}^2
	\bigr)^2\\
	  &=
	  \sum_{\vect{h} \in \generatingSet(\mat{M}^\tT)}
	  \sum_{\vect{z} \in \Z^d}
	  m^2
	  \frac{u_{\Stiffness^0}}{l_{\Stiffness^0}}
	  \norm{\hat G_\vect{h} c_{\vect{h} + \mat{M}^\tT \vect{z}}(f)}^2
	  m^{-2}.
	\end{align*}
	With another application of the Parseval equation the desired
	estimate
	\begin{equation*}
    \begin{split}
  	  \norm[\bigl]{\Gamma^p \Stiffness^0 &\tensorProd \gamma}[\Sym_d(L^2(\T^d))]^2
	  \\&\leq
	  \sum_{\vect{h} \in \generatingSet(\mat{M}^\tT)}
	  \sum_{\vect{z} \in \Z^d}
	  \frac{u_{\Stiffness^0}}{l_{\Stiffness^0}}
	  \norm{\hat G_\vect{h} c_{\vect{h} + \mat{M}^\tT \vect{z}}(f)}^2
	  =
	  \frac{u_{\Stiffness^0}}{l_{\Stiffness^0}}
	  \norm{\gamma}[\Sym_d(L^2(\T^d))]^2.
    \end{split}
	\end{equation*}
  is obtained.
  \end{proof}
\end{theorem}

The computation of the Green operator $\Gamma^p$ on $V^f_\mat{M}$ involves
computing the value of the series 
\begin{equation*}
  m\sum_{\vect{z} \in \Z^d}
  \hat{\Gamma}^0_{\vect{h} + \mat{M}^\tT \vect{z}}
  \abs{c_{\vect{h} + \mat{M}^\tT \vect{z}}(f)}^2
\end{equation*}
for all $\vect{h} \in \generatingSet(\mat{M}^\tT)$. This evaluation simplifies
for functions $f$ having compact support in the frequency domain and where the
series reduces to a sum over finitely many terms. In this case the operator can
be evaluated exactly, i.e. without introducing any additional numerical error.
This is the case for example for de la Vall\'ee Poussin means, where each sum only consists of up to \(4\) terms.

Functions that have compact support in space no longer allow for an exact
evaluation of $\Gamma^p$. An example are Box splines in space domain which can
be interpreted as finite elements integrated by only one quadrature point. In
addition, the Box splines allow for finite elements which have different degrees
of differentiability in directions other than the grid.

Brisard and Dormieux~\cite{BrisardDormieux2010Framework} derive a Green
operator from an energy based formulation using a discretisation with
element-wise constant finite elements. Their Green operator
corresponds to using a Box spline of order zero in the approach here
and is thus contained in the framework.

\subsection{Discretisation of the Lippmann-Schwinger equation}

With the definition of the periodised Green operator $\Gamma^p$ we can now
proceed to derive a corresponding discretisation of the Lippmann-Schwinger
equation. This derivation is split into two theorems.
The discretised version of the space $\Curlfree\bigl(\T^d\bigr)$ is given by
\begin{equation*}
  \Curlfree_\mat{M}^f\bigl( \T^d \bigr) \coloneqq
  \bigl\lbrace
  \varepsilon_\mat{M} \in \Sym_d\bigl(V_\mat{M}^f \bigr): \varepsilon_\mat{M} \text{
  interpolates } \varepsilon \in \Curlfree\bigl(\T^d\bigr) \text{ on }
  \Pattern(\mat{M})
  \bigr\rbrace.
\end{equation*}
This space allows to state the discretised version of the PDE
correctly. For the following theorems we assume that $\bigl(\Stiffness -
\Stiffness^0 \bigr) \tensorProd (\varepsilon_\mat{M} + \varepsilon^0) \in
\Sym_d\bigl(\Wiener(\T^d)\bigr)$ for $\varepsilon_\mat{M} \in
\Curlfree_\mat{M}^f\bigl(\T^d\bigr)$ so an interpolation on $\Sym_d\bigl(
V_\mat{M}^f \bigr)$ is possible with
\begin{equation*}
  \Bigl(
    \bigl(\Stiffness - \Stiffness^0 
    \bigr)\tensorProd \bigl(\varepsilon + \varepsilon^0\bigr)
  \Bigr)(\vect{x}) =
  \sum_{\vect{y} \in \Pattern(\mat{M})}
  B_\vect{y} \Translate(\vect{y}) \Fundamental_\mat{M}(\vect{x})
\end{equation*}
for all $\vect{x} \in \Pattern(\mat{M})$.

Additionally, a test functions $\gamma$ can be written as
\begin{equation}
  \label{eq:testfunction_factors}
  \gamma = \sum_{\vect{y} \in \Pattern(\mat{M})} 
  G_\vect{y} \Translate(\vect{y}) f
\end{equation}
with \(\vect{G}
= (G_{\vect{y}})_{\vect{y}\in\Pattern(\mat{M})}\in\mathbb C^m\)
and its discrete Fourier transform by~\(\vect{\hat G} = (\hat G_{\vect{h}})_{\vect{h}\in\generatingSet(\mat{M}^\tT)} = \Fourier(\mat{M})\vect{G}\).

Let the strain  $\varepsilon_\mat{M} \in \Curlfree_\mat{M}^f\bigl(\T^d\bigr)$ be
written in terms of translates of the fundamental interpolant as
\begin{equation*}
  \varepsilon_\mat{M} =
  \sum_{\vect{y} \in \Pattern(\mat{M})}
  E_\vect{y}
  \Translate(\vect{y}) \Fundamental_\mat{M}
\end{equation*}
and again the discrete Fourier transform of the coefficient vector
as~\(\vect{\hat E} = \Fourier(\mat{M})\vect{E}\).

\begin{theorem}
  \label{thm:pde_interpolation}
  Let the translates of $f$ be orthonormal, let $\bigl(\Stiffness -
  \Stiffness^0 \bigr) \tensorProd (\varepsilon_\mat{M} + \varepsilon^0) \in
  \Sym_d\bigl(\Wiener(\T^d)\bigr)$, let $\gamma \in
  \Sym_d\bigl(V^f_\mat{M} \bigr)$, and $\varepsilon_\mat{M} \in
  \Curlfree_\mat{M}^f\bigl(\T^d\bigr)$. Then it holds
  \begin{equation}\label{eq:pde_interpolation}
    \begin{split}
          \bigl\langle
    \varepsilon_\mat{M} &+
    \Gamma^0 \bigl(
      \Stiffness - \Stiffness^0
    \bigr) \tensorProd
    \bigl(
      \varepsilon_\mat{M} + \varepsilon^0
    \bigr),
    \gamma
    \bigr\rangle
    \\&=
      \sum_{\vect{h} \in \generatingSet(\mat{M}^\tT)}
      \frac{1}{m}
      \bigl\langle
      \hat E_\vect{h} \hat{a}_\vect{h},
	\hat G_\vect{h}
	\bigr\rangle + 
      \Biggl\langle
       \hat{a}_\vect{h}
      \Bigl[
	\bigl\lbrace
	\hat{\Gamma}^0_\vect{k} \abs{c_\vect{k}(f)}^2
	\bigr\rbrace_{\vect{k} \in \Z^d}
      \Bigr]^\mat{M}_\vect{h}
      \hat B_\vect{h},
      \hat G_\vect{h}
      \Biggr\rangle.
    \end{split}
  \end{equation}
  \begin{proof}
    Starting with the left-hand side of~\eqref{eq:pde_interpolation} applying
    the Parseval equation~\eqref{eq:parseval} to transform it to Fourier space
    yields the equal form
    \begin{equation*}
      \sum_{\vect{k} \in \Z^d}
      \bigl\langle
      c_\vect{k}(\varepsilon_\mat{M}) +
      \hat{\Gamma}^0_\vect{k} 
      c_\vect{k} \Bigl(
	\bigl(\Stiffness - \Stiffness^0 \bigr) \tensorProd
	 \bigl( \varepsilon_\mat{M} + \varepsilon^0 \bigr)
      \Bigr),
      c_\vect{k}(\gamma)
      \bigr\rangle,
    \end{equation*}
    where we make use of~\eqref{eq:gamma}.
    Equation~\eqref{eq:inTranslates:ck}
    together with the formula for the translate coefficients of the fundamental
    interpolant~\eqref{eq:fundamental_coeff} and the splitting $\vect{k} =
    \vect{h} + \mat{M}^\tT \vect{z}$ with $\vect{k},\vect{z} \in \Z^d$ and
    $\vect{h} \in \generatingSet(\mat{M}^\tT)$ results in the expressions
    \begin{align*}
      c_\vect{k} 
      \Bigl(
	\bigl(\Stiffness - \Stiffness^0 \bigr) \tensorProd
	\bigl(\varepsilon_\mat{M} + \varepsilon^0\bigr)
      \Bigr) &= 
      \hat B_\vect{h} \hat{a}_\vect{h}
      c_{\vect{h} + \mat{M}^\tT \vect{z}}(f),\\
      c_\vect{k}(\gamma) &= 
      \hat G_\vect{h} c_{\vect{h} + \mat{M}^\tT \vect{z}}(f),\\
      c_\vect{k}(\varepsilon_\mat{M}) &= 
      \hat E_\vect{h} \hat{a}_\vect{h} c_{\vect{h} + \mat{M}^\tT \vect{z}}(f).
    \end{align*}
    Inserting these equations into the expression above yields
    \begin{equation*}
      \begin{split}
      \sum_{\vect{h} \in \generatingSet(\mat{M}^\tT)}
      \sum_{\vect{z} \in \Z^d}
      \bigl\langle&
      \hat E_\vect{h} \hat{a}_\vect{h} c_{\vect{h} + \mat{M}^\tT \vect{z}}(f),
      \hat G_\vect{h} c_{\vect{h} + \mat{M}^\tT \vect{z}}(f)
      \bigr\rangle
      \\&\quad+ 
      \bigl\langle
      \Gamma^0_{\vect{h} + \mat{M}^\tT \vect{z}}
      \hat B_\vect{h} \hat{a}_\vect{h}
      c_{\vect{h} + \mat{M}^\tT\vect{z}}(f),
      \hat G_\vect{h} c_{\vect{h} + \mat{M}^\tT\vect{z}}(f)
      \bigr\rangle.
      \end{split}
    \end{equation*}
    Collecting the terms depending on $\vect{z}$ and employing the bracket
    sums~\eqref{eq:bracketsum} to simplify the expression, one obtains
    \begin{equation*}
      \sum_{\vect{h} \in \generatingSet(\mat{M}^\tT)}
      \Bigl\langle
      \hat E_\vect{h} \hat{a}_\vect{h}
	\bigl[
	  \bigl\lbrace
	  \abs{ c_\vect{k}(f) }^2
	  \bigr\rbrace
	\bigr]_\vect{h}^\mat{M},
	\hat G_\vect{h}
	\Bigr\rangle + 
      \Biggl\langle
       \hat{a}_\vect{h}
      \Bigl[
	\bigl\lbrace
	\hat{\Gamma}^0_\vect{k} \abs{c_\vect{k}(f)}^2
	\bigr\rbrace_{\vect{k} \in \Z^d}
      \Bigr]^\mat{M}_\vect{h}
      \hat B_\vect{h},
      \hat G_\vect{h}
      \Biggr\rangle.
    \end{equation*}
    Since by assumption translates of $f$ are orthonormal they fulfil
    $
	\bigl[
	  \bigl\lbrace
	  \abs{ c_\vect{k}(f) }^2
	  \bigr\rbrace
	\bigr]_\vect{h}^\mat{M} = \frac{1}{m}
	$ by Lemma~\ref{lem:FI:Props}~\ref{lem:FI:Props:ONB} and we get
  the desired result.
  \end{proof}
\end{theorem}

The following theorem states the result of a Galerkin projection
of~\eqref{eq:weak_pde_projected} onto the space of translates.

\begin{theorem}
  \label{thm:pde_discretization}
  Let the translates of $f$ be orthonormal, let $\Stiffness \in
  \SSym_d\bigl(A\bigl(\T^d\bigr)\bigr)$, and let $ \varepsilon_\mat{M} \in
  \Curlfree_\mat{M}^f\bigl(\T^d\bigr)$. Then $\varepsilon_\mat{M}$ fulfils the
  weak form
  \begin{equation}
    \label{eq:weak_pde_projected_1}
    \bigl\langle
    \varepsilon_\mat{M} +
    \Gamma^0 \bigl(
      \Stiffness - \Stiffness^0
    \bigr) \tensorProd
    \bigl(
      \varepsilon_\mat{M} + \varepsilon^0
    \bigr),
    \gamma
    \bigr\rangle = 0
  \end{equation}
  for all $\gamma \in \Sym_d\bigl( V_\mat{M}^f \bigr)$ if and only
  if
  \begin{equation}
    \label{eq:pde_coefficients}
    \sum_{\vect{y} \in \Pattern(\mat{M})}
    \Bigl(
      E_\vect{y} + 
      \Gamma^\p_\vect{y}
      \bigl(\Stiffness(\vect{y}) - \Stiffness^0 \bigr) \tensorProd
      \bigl( E_\vect{y} + \varepsilon^0 \bigr)
    \Bigr)
    \bigl(\Translate(\vect{y}) \Fundamental_\mat{M} \bigr)(\vect{x}) = \vect{0}
  \end{equation}
  for all $\vect{x} \in \Pattern(\mat{M})$, where $\Gamma^\p_\vect{y}$ is
  defined in Definition~\ref{def:gamma_tilde_fourier}.
  \begin{proof}
    With $\varepsilon_\mat{M} \in  \Curlfree_\mat{M}^f\bigl(\T^d\bigr) \subset
    \Sym_d\bigl( V_\mat{M}^f \bigr) \subset \Sym_d\bigl(A(\T^d)\bigr)$ and
    $\Stiffness \in \SSym_d\bigl(A\bigl(\T^d\bigr)\bigr)$ it follows that
    $\bigl( \Stiffness - \Stiffness^0 \bigr) \tensorProd (\varepsilon_\mat{M} +
    \varepsilon^0) \in \Sym_d\bigl(A(\T^d)\bigr)$ and hence the assumptions of
    Theorem~\ref{thm:pde_interpolation} are fulfilled.
    Therefore~\eqref{eq:weak_pde_projected_1} is equivalent to
    \begin{equation}
      \label{eq:pde_interpolation_eq}
      \sum_{\vect{h} \in \generatingSet(\mat{M}^\tT)}
      \frac{1}{m}
      \bigl\langle
      \hat E_\vect{h} \hat{a}_\vect{h},
	\hat G_\vect{h}
	\bigr\rangle + 
      \Biggl\langle
       \hat{a}_\vect{h}
      \Bigl[
	\bigl\lbrace
	\hat{\Gamma}^0_\vect{k} \abs{c_\vect{k}(f)}^2
	\bigr\rbrace_{\vect{k} \in \Z^d}
      \Bigr]^\mat{M}_\vect{h}
      \hat B_\vect{h},
      \hat G_\vect{h}
      \Biggr\rangle = 0,
    \end{equation}
    with the notation from~\eqref{eq:testfunction_factors}.
    A necessary and sufficient condition for~\eqref{eq:pde_interpolation_eq} to
    hold true is that it is fulfilled for all
    \begin{equation*}
      \hat G_{\vect{h},\vect{y},p,q} \coloneqq \beta_{pq} \bigl(
	\alpha_p \alpha_q^\tT +
      \alpha_q \alpha_p^\tT \bigr) \e^{-2 \pi \im \vect{h}^\tT \vect{y}}
    \end{equation*}
    for all $\vect{h} \in \generatingSet(\mat{M}^\tT)$ and $\vect{y} \in
    \Pattern(\mat{M})$ and $p,q \in \lbrace 1,\dots,d \rbrace$.
    The vector
    $\alpha_p \in \R^d$ denotes the $p$-th unit vector and
    $\beta_{pq} \coloneqq 1-\frac12 \delta_{pq}$ normalizes the resulting
    matrix. This parametrization is the trigonometric basis of
    $\Sym_d\bigl(V_\mat{M}^f\bigr)$
    on the pattern $\Pattern(\mat{M})$.

    Hence an equivalent condition stems from looking
    at~\eqref{eq:pde_interpolation_eq}
    component-wise, i.e.
    \begin{equation*}
      \sum_{\vect{h} \in \generatingSet(\mat{M}^\tT)}
      \frac{1}{m} \hat E_\vect{h} \hat{a}_\vect{h} e^{2 \pi \im \vect{h}^\tT \vect{y}} + 
      \sum_{\vect{h} \in \generatingSet(\mat{M}^\tT)}
      \hat{a}_\vect{h}
      \Bigl[
	\bigl\lbrace
	\hat{\Gamma}^0_\vect{k} \abs{c_\vect{k}(f)}^2
	\bigr\rbrace_{\vect{k} \in \Z^d}
      \Bigr]^\mat{M}_\vect{h}
      \hat B_\vect{h}
      \e^{2 \pi \im \vect{h}^\tT \vect{y}} = \vect{0},
    \end{equation*}
 bearing in mind the necessary complex conjugate, for all $\vect{y} \in
 \Pattern(\mat{M})$.
 This, however, is an inverse discrete Fourier transform on the pattern
 $\Pattern(\mat{M})$ and together with
 Definition~\ref{def:gamma_tilde_fourier} and setting
 \begin{equation*}
   \hat{\tilde B}_\vect{h} \coloneqq 
   m
   \Bigl[
     \bigl\lbrace
     \hat{\Gamma}^0_\vect{k} \abs{c_\vect{k}(f)}^2
     \bigr\rbrace_{\vect{k} \in \Z^d}
   \Bigr]^\mat{M}_\vect{h}
   \hat B_\vect{h}
   = \hat{\Gamma}^\p_\vect{h} \hat B_\vect{h}
 \end{equation*}
 yields
    \begin{align*}
      \frac{1}{m}&\sum_{\vect{h} \in \generatingSet(\mat{M}^\tT)}
      \bigl(
	\hat E_\vect{h} +
	\hat{\tilde B}_\vect{h} 
      \bigr)
      \hat{a}_\vect{h} \e^{2 \pi \im \vect{h}^\tT \vect{y}}
      \\&=
      \frac{1}{m}\sum_{\vect{h} \in \generatingSet(\mat{M}^\tT)}
      \biggl(
	\hat E_\vect{h} \hat{a}_\vect{h}+ 
	\hat{a}_\vect{h} \Stiffness^0 m
	\Bigl[
	  \bigl\lbrace
	  \hat{\Gamma}^0_\vect{k} \abs{c_\vect{k}(f)}^2
	  \bigr\rbrace_{\vect{k} \in \Z^d}
	\Bigr]^\mat{M}_\vect{h} \tensorProd
	\hat B_\vect{h}
      \biggr)
      \e^{2 \pi \im \vect{h}^\tT \vect{y}}\\
      &=
      \frac{1}{m}\sum_{\vect{h} \in \generatingSet(\mat{M}^\tT)}
      \Bigl(
	\hat E_\vect{h} \hat{a}_\vect{h} +
	\hat{a}_\vect{h}
	\hat{\Gamma}^\p_\vect{h}
	\hat B_\vect{h}
      \Bigr)
      \e^{2 \pi \im \vect{h}^\tT \vect{y}}.
    \end{align*}
    The coefficients
    $\bigl(
      \tilde B_\vect{y}
      \bigr)_{\vect{y} \in \Pattern(\mat{M})}
      \coloneqq \frac{1}{\sqrt{m}} \overline{\Fourier}(\mat{M})^\tT
    \bigl(
    \hat{\tilde B}_\vect{h}
    \bigr)_{\vect{h} \in \generatingSet(\mat{M}^\tT)}$
    can now be interpreted as coefficients of
    translates of the fundamental interpolant, i.e. it holds
 \begin{equation*}
   \sum_{\vect{y} \in \Pattern(\mat{M})}
   \bigl(
     E_\vect{y} + 
     \tilde B_\vect{y}
   \bigr)
   \bigl(\Translate(\vect{y}) \Fundamental_\mat{M}\bigr)(\vect{x}) =
   \vect{0}
 \end{equation*}
 for all~$\vect{x} \in \T^d$.
 By Definition~\ref{def:gamma_tilde_fourier} the operator $\Gamma^\p$ acts as a
 Fourier multiplier with Fourier coefficients~\eqref{eq:gamma_tilde_fourier}.
 This transforms the above equation to
 \begin{equation*}
   \sum_{\vect{y} \in \Pattern(\mat{M})}
   \Bigl(
     E_\vect{y} + 
     \Gamma^\p_\vect{y} B_\vect{y}
   \Bigr)
   \bigl(\Translate(\vect{y}) \Fundamental_\mat{M}\bigr)(\vect{x}) = \vect{0}.
 \end{equation*}
 The coefficients $B_\vect{y}$ were chosen such that they coincide with
 the function values of $
 \bigl(
   \Stiffness(\vect{y}) - \Stiffness^0
 \bigr)
 \tensorProd
 \bigl(
 \varepsilon_\mat{M}(\vect{y}) + \varepsilon^0
 \bigr)$
 at points $\vect{y} \in \Pattern(\mat{M})$. Likewise, $\varepsilon_\mat{M}$
 coincides in the points $\vect{y} \in \Pattern(\mat{M})$ with the coefficients
 $E_\vect{y}$. Inserting these relations one obtains
 \begin{equation*}
   \sum_{\vect{y} \in \Pattern(\mat{M})}
   \Bigl(
     E_\vect{y} + 
     \Gamma^\p_\vect{y}
     \bigl(\Stiffness(\vect{y}) - \Stiffness^0 \bigr) \tensorProd
     \bigl( E_\vect{y} + \varepsilon^0 \bigr)
   \Bigr)
   \bigl(\Translate(\vect{y}) \Fundamental_\mat{M} \bigr)(\vect{x}) = \vect{0}
 \end{equation*}
 which yields the desired result.
  \end{proof}
\end{theorem}

When discretising the PDE~\eqref{eq:weak_pde_projected} in a similar way, one
arrives at the following discretised form.
\begin{theorem}
  \label{thm:ls_discretization}
  Let the translates of $f$ be orthonormal and let $\Stiffness \in
  \SSym_d\bigl(A\bigl(\T^d\bigr)\bigr)$ and let $ \varepsilon_\mat{M}
  \in \Curlfree_\mat{M}^f\bigl(\T^d\bigr)$.
  Then $\varepsilon_\mat{M}$ fulfils the weak form~\eqref{eq:weak_pde_projected}
  \begin{equation*}
    \bigl\langle
    \Stiffness^0
    \Gamma^0
      \Stiffness
    \tensorProd
    \bigl(
      \varepsilon_\mat{M} + \varepsilon^0
    \bigr),
    \gamma
    \bigr\rangle = 0
  \end{equation*}
for all $\gamma \in \Sym_d\bigl( V_\mat{M}^f \bigr)$ if and only if
  \begin{equation}\label{eq:pde_coefficients_variational}
    \sum_{\vect{y} \in \Pattern(\mat{M})}
      \Stiffness^0 \Gamma^\p_\vect{y}
      \Stiffness(\vect{y}) \tensorProd
      \bigl( E_\vect{y} + \varepsilon^0 \bigr)
      \bigl(\Translate(\vect{y}) \Fundamental_\mat{M} \bigr)(\vect{x}) = \vect{0}
  \end{equation}
  for all $\vect{x} \in \Pattern(\mat{M})$.
  \begin{proof}
    The proof follows the same steps as for the
    Theorems~\ref{thm:pde_interpolation} and~\ref{thm:pde_discretization}
    and we omit it therefore.
  \end{proof}
\end{theorem}

In Remark~\ref{rem:Gamma_identity} we already mentioned that in the continuous
case the Lippmann-Schwinger equation~\eqref{eq:lippmann-schwinger} and the
variational equation~\eqref{eq:weak_pde_projected} coincide. This is, as
shown in~\cite{Vondrejc2014}, also the case when using trigonometric collocation
for the discretisation. With the equations~\eqref{eq:pde_coefficients}
and~\eqref{eq:pde_coefficients_variational} using spaces of translates this is
in general no longer the case. When looking at the
identity~\eqref{eq:projection_identity} one can see this rather quickly.

\begin{remark}
  \label{rem:constants}
  For $\varepsilon_\mat{M} \in \Curlfree_\mat{M}^f\bigl(\T^d\bigr)$ it holds true
  that
  \begin{equation*}
    \varepsilon_\mat{M} - \Gamma^\p \Stiffness^0 \tensorProd \varepsilon_\mat{M} =
    \varepsilon^0
  \end{equation*}
  almost everywhere if and only if \(V_{\mat{M}}^f=V_{\mat{M}}^{D_{\mat{M}}}\),
  i.e.~if $f$ can be written as a sum of translates of the Dirichlet kernel $D_{\mat{M}}$.
  \begin{proof}
    The proof in~\cite[Proposition 3]{Vondrejc2014} uses that $\Gamma^0
    \Stiffness^0$ projects a constant function $g$ onto the function that is $0$
    almost everywhere, i.e. it is only characterized by the Fourier coefficient
    $c_\vect{0}(g)$. When interpolating in $V^f_\mat{M}$ this is in general no
    longer the case and an application of $\Gamma^\p \Stiffness^0$ does not
    result in the zero function.
  \end{proof}
\end{remark}

\begin{figure}\centering
  \includegraphics{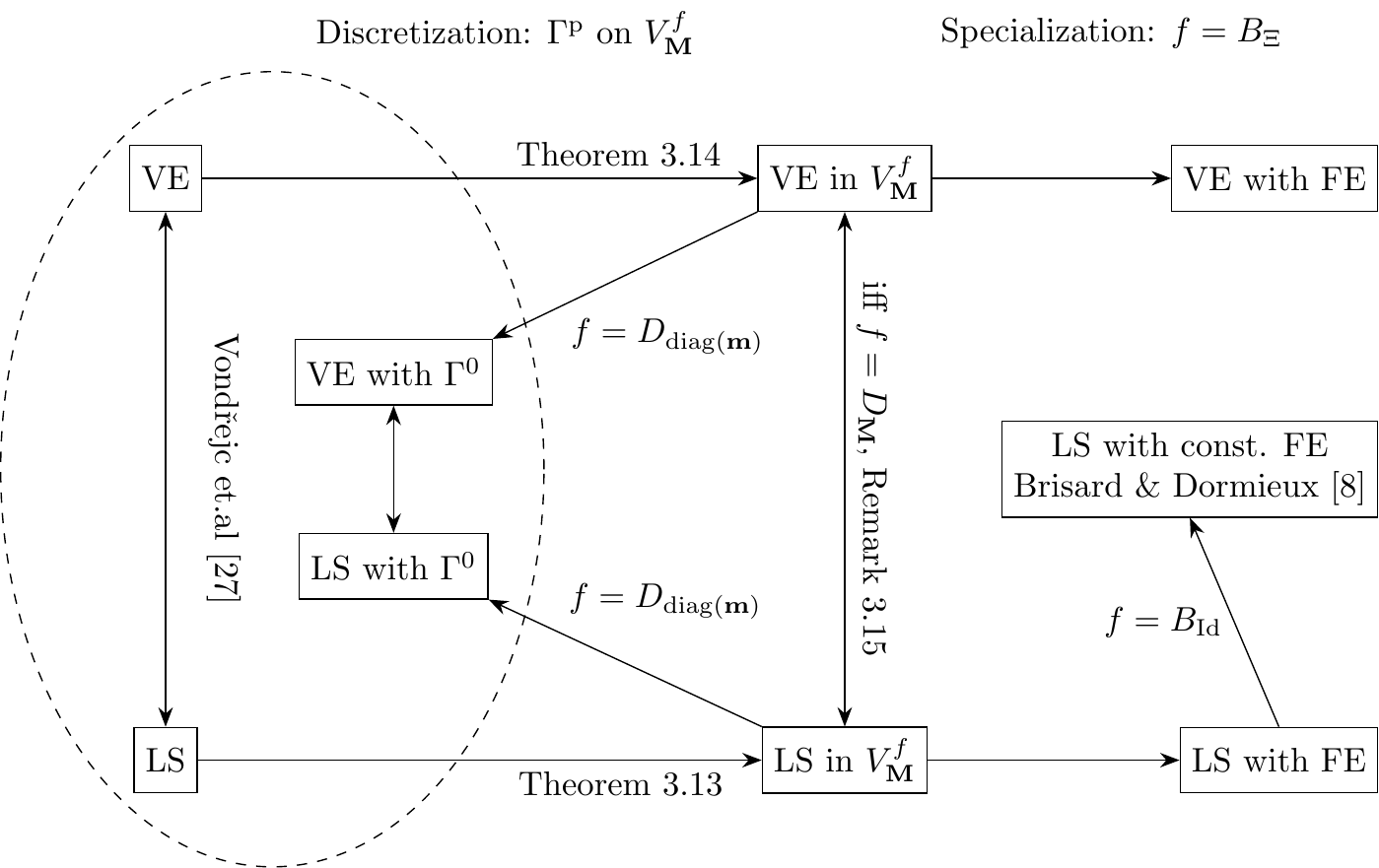}
  \caption[]{A diagram of connections between the Lippmann-Schwinger equation
    (LS) in~\eqref{eq:lippmann-schwinger} and the variational equation (VE)
  in~\eqref{eq:weak_pde_projected} for different discretisations. The term 
  $\diag(\vect{m})$ with $\vect{m} \in \mathbb N^d$ denotes a diagonal matrix and
thus $D_{\diag(\vect{m})}$ is the Dirichlet kernel on a tensor product grid.}
  \label{fig:diagram}
\end{figure}

The connections between the variational
formulation~\eqref{eq:weak_pde_projected} and the Lippmann-Schwinger
equation~\eqref{eq:lippmann-schwinger} and their discretisations on spaces of
translates are summarized in Figure~\ref{fig:diagram}.
The continuous equations are shown to be equivalent in~\cite[Proposition
3]{Vondrejc2014}. The same holds also true for the equations discretised on
$V_\mat{M}^f$ if and only if $V_{\mat{M}}^f=V_{\mat{M}}^{D_{\mat{M}}}$
For the special case~$f = D_\mat{M}$ and $\mat{M}$ a diagonal matrix, i.e. for
a tensor product grid,
this equivalence is already proven in~\cite[Proposition 12]{Vondrejc2014}. Box
splines allow to solve the equations in terms of (simplified) finite elements,
which generalizes the constant finite element approach
of~\cite{BrisardDormieux2010Framework}. This emerges when discretising the
Lippmann-Schwinger equation with $f=B_{\Id}$.

\section{Numerics}
\label{sec:numerics}

In this section we study the effect of choosing different functions for
the translation invariant spaces to illustrate the capabilities of the
generalization presented in this paper. In
the publication~\cite{BergmannMerkert2016} the authors study the
influence of different patterns on the solution quality and their numerical
effects.

\begin{figure}\centering
  \begin{subfigure}[c]{.49\textwidth}\centering
    \includegraphics[width=.98\textwidth]{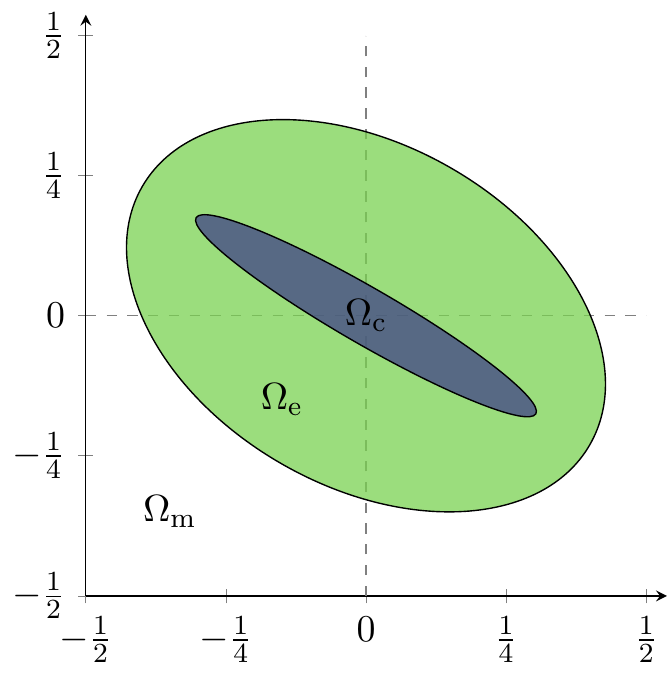} 
  \end{subfigure}
  \begin{subfigure}[c]{.49\textwidth}\centering
    \includegraphics[width=.98\textwidth]{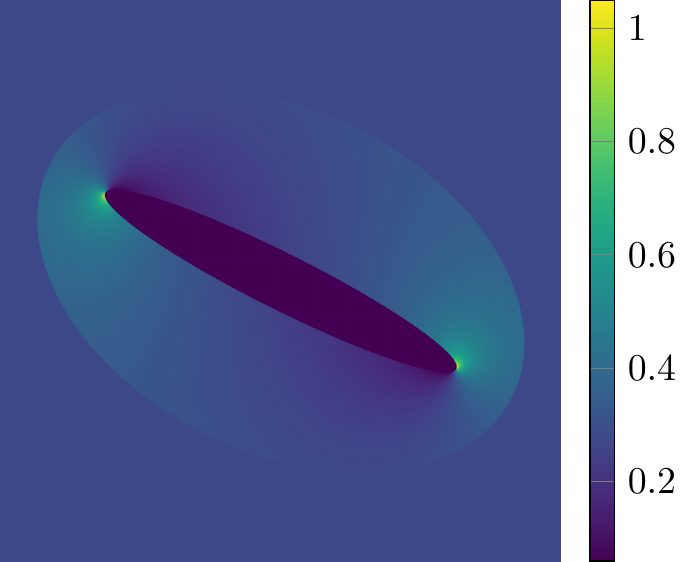} 
  \end{subfigure}
  \caption{A schematic of the generalized Hashin structure (left) and the
    analytic solution for the first component of the strain $\varepsilon_{11}$
  (right).}
  \label{fig:hashin}
\end{figure}

A prototypical structure that is introduced in~\cite{BergmannMerkert2016} is
the generalized Hashin structure, a geometry that is based on publications of
Milton, see~\cite{Milton2002}. It consists of two confocal ellipses embedded in
a surrounding material, see Figure~\ref{fig:hashin}. The centre ($\Omega_c$) and
coating ($\Omega_e$) ellipses have isotropic behaviour and the matrix material
($\Omega_m$) is built in such a way that it is unaffected by the inclusion for a
chosen macroscopic strain~$\varepsilon^0$. For this special kind of structure
analytic expressions for the strain field $\varepsilon$ and the action of the
effective matrix $\Stiffness^\eff \tensorProd \varepsilon^0$ are known.

In the following we take exactly the same parameters as
in~\cite{BergmannMerkert2016}, i.e. for the ellipses with parametrization
\begin{equation*}
  \frac{x_1^2}{c_1^2+\rho}
  + \frac{x_2^2}{c_2^2+\rho} = 1
\end{equation*}
and $(x_1, x_2)^\tT \in \T^2$ we choose $c_1 = 0.05$, $c_2 = 0.35$ and for the
inner and outer ellipsis $\rho = 0$ and $\rho = 0.09$, respectively. The
structure is then rotated by $60^\circ$. For the inner and outer ellipsis we
choose isotropic material laws with Poisson's ratio $\nu = 0.3$ in both ellipses
and the matrix material and Young's moduli $E = 1$ and $E = 10$ for the inner
and outer ellipsis, respectively. The material law for the surrounding matrix
material can then be determined by the formulae in~\cite[Section
4.2]{BergmannMerkert2016}.

A solution for the first component of the strain field $\varepsilon_{11}$ is
depicted in Figure~\ref{fig:hashin} and taken for instructional purposes
directly from~\cite{BergmannMerkert2016}.

All numerical results in this section were obtained by
solving~\eqref{eq:pde_coefficients} with a fixed-point iteration on
$E_\vect{y}$ as described in~\cite{MoulinecSuquet1998} up to an relative error
of $10^{-10}$ using a Cauchy criterion.

\subsection{De la Vallée Poussin means}

In~\cite{BergmannMerkert2016} the authors study the influence of the pattern
matrix $\mat{M}$ on the solution field and the quality of the effective matrix
$\Stiffness^\eff$. We take the following pattern matrices:
\begin{equation*}
    \mat{M}_1 = 
    \begin{pmatrix} 128 & 0\\0 & 128 \end{pmatrix},\quad
    \mat{M}_2 = 
    \begin{pmatrix} 64 & 64\\-64 & 64 \end{pmatrix},\quad
    \mat{M}_3 = 
    \begin{pmatrix} 128 & 272\\0 & 128 \end{pmatrix}.
  \end{equation*}
  They correspond to a tensor product grid ($\mat{M}_1$), a tensor
  product grid rotated by $45^\circ$ ($\mat{M}_2$) and the minimal $\ell^2$-error
  ($\mat{M}_3$) achieved. The matrices $\mat{M}_1$ and $\mat{M}_3$ have a
  determinant of $2^{14}$ whereas the so-called quincux pattern
  $\Pattern(\mat{M}_2)$ has $2^{13}$ sampling points.

  First we study how using de la Vall\'ee Poussin means changes the quality of
  the effective stiffness matrix and the strain field. The Box splines used in
  frequency domain to define the means allow to use different slopes in each
  direction of the de la Vall\'ee Poussin mean and thus one can expect to reduce
  the Gibbs phenomenon in different directions. For the following study the
  functions $f_{\mat{M}_i,\alpha}$ are parametrized with $\alpha = (\alpha_1,
  \alpha_2)^\tT$ and $\alpha_1, \alpha_2 \in [0,0.5]$ for $i = 1,\dots,3$.

  The parameters $\alpha_1$ and $\alpha_2$ correspond to damping the Fourier
  coefficients of the de la Vall\'ee Poussin mean along the directions
  $\mat{M}^\tT ( 1, 0)^\tT$ and $\mat{M}^\tT ( 0,1)^\tT$, respectively. In the
  space domain they introduce a better localization~\cite{GohGoodman2004} along
  $\mat{M}^{-1} (0, 1)^\tT$ and $\mat{M}^{-1} ( 1,0)^\tT$, respectively.

\begin{figure}
  \begin{subfigure}[b]{.48\textwidth}
    \centering
    \includegraphics[width=.98\textwidth]{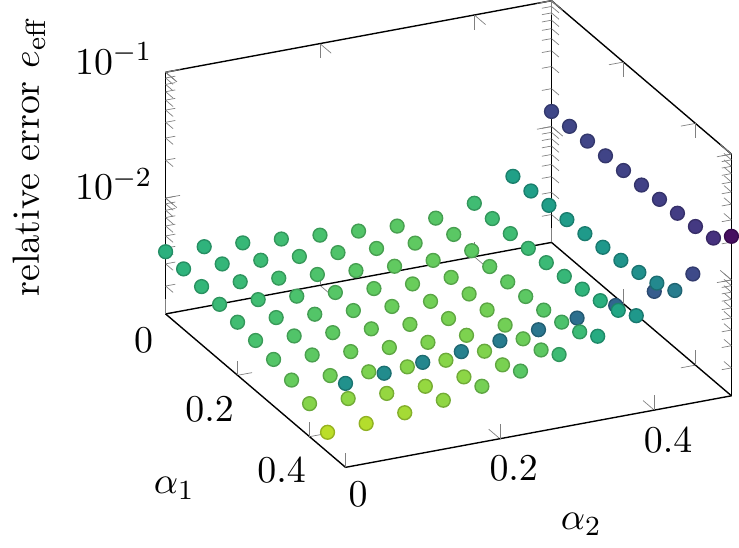}
  \end{subfigure}
  \hfill
  \begin{subfigure}[b]{.48\textwidth}
    \centering
    \includegraphics[width=.98\textwidth]{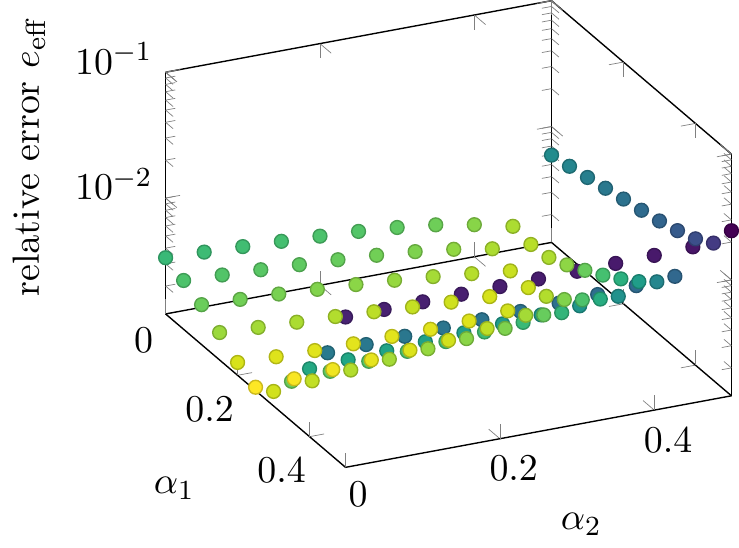}
  \end{subfigure}\\
  \centering
  \begin{subfigure}[b]{.48\textwidth}
    \centering
    \includegraphics[width=.98\textwidth]{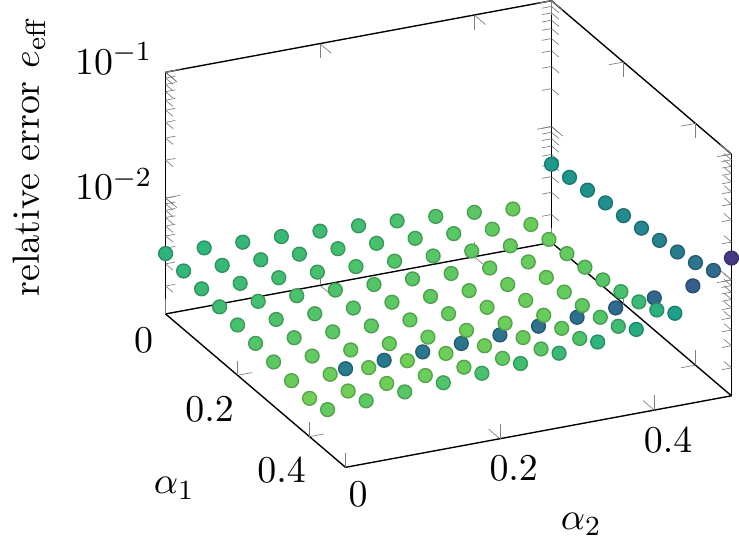}
  \end{subfigure}
  \caption[]{
    Relative error $e_\eff$ depending on the slopes of the De la Vallée Poussin
    means with parameters $\alpha_1$ and $\alpha_2$ in the directions of the
    unit cell. The matrices inducing the pattern are from left to right, top to
    bottom:
    $\mat{M}_1 = 
    \bigl(\begin{smallmatrix} 128 & 0\\0 & 128 \end{smallmatrix}\bigl)$,
    $\mat{M}_2 = 
    \bigl(\begin{smallmatrix} 64 & 64\\-64 & 64 \end{smallmatrix}\bigl)$,
    $\mat{M}_3 = 
    \bigl(\begin{smallmatrix} 128 & 272\\0 & 128 \end{smallmatrix}\bigl)$.
  }
  \label{fig:dlvp_effective}
\end{figure}

  The result of these experiments regarding the effective stiffness matrix is
  shown in Figure~\ref{fig:dlvp_effective}. 
  The relative effective error
  \begin{equation*}
  e_\eff \coloneqq
  \Biggl\lVert 
  \Stiffness^\eff \tensorProd \varepsilon^0 - 
  \sum_{\vect{y} \in \Pattern(\mat{M})}
  \Stiffness(\vect{y}) \tensorProd E_\vect{y}
  \Biggr\rVert
  \bigl\lVert
  \Stiffness^\eff \tensorProd \varepsilon^0
  \bigr\rVert^{-1}
  \end{equation*}
  is depicted with $E_\vect{y}$ from~\eqref{eq:pde_coefficients}.
  Parameters $\alpha_1 = \alpha_2 = 0$ correspond to the modified Dirichlet
  kernel $f_{\mat{M},\vect{0}}$ and $\alpha_1 = \alpha_2 = 0.5$ correspond to
  the Fej\'er kernel $f_{\mat{M},\frac12 \vect{1}}$.

  For $\mat{M}_1$, the error in the effective stiffness matrix changes from
  $0.0037$ for the modified Dirichlet kernel to the optimum at $\alpha_1 = 0.45$
  and $\alpha_2 = 0$ with $e_\eff = 0.0017$. For $\mat{M}_2$ it changes from
  $0.0034$ to $e_\eff = 0.0013$ for $\alpha_1 = 0.25$ and $\alpha_2 = 0$. The
  solution for $\mat{M}_3$ can be improved upon by de la Vall\'ee Poussin means
  with $\alpha_1 = 0.4$ and $\alpha_2 = 0$ which results in $e_\eff = 0.0036$ of
  the Fourier approach on patterns being reduced to $e_\eff = 0.0024$. When
  approaching either $\alpha_1 = 0.5$ or $\alpha_2 = 0.5$ the value of $e_\eff$
  drastically increases.
  The errors for the solution using the Dirichlet kernel and the modified
  Dirichlet kernel, i.e. for $\alpha_1 = \alpha_2 = 0$, are almost the same and
  thus not marked in the plots above.

  The study in~\cite{BergmannMerkert2016} suggests that by changing the pattern
  one can improve the quality of the solution and the effective stiffness
  matrix. The experiments from above show that one can improve these results
  even further by extending the theory to translation invariant spaces.
  Especially for the tensor product grid, i.e. for data given on a regular voxel
  grid like for example from a computer tomography image, a suitable choice of
  the space of translates could here reduce the error in the effective stiffness
  by $55\%$.

\begin{figure}
  \begin{subfigure}{\textwidth}\centering
    \includegraphics[scale=0.8]{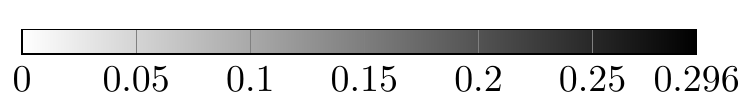}
  \end{subfigure}
    \\[.25\baselineskip]
  \begin{subfigure}[t]{.32\textwidth}\centering
    \includegraphics[width=0.95\textwidth]{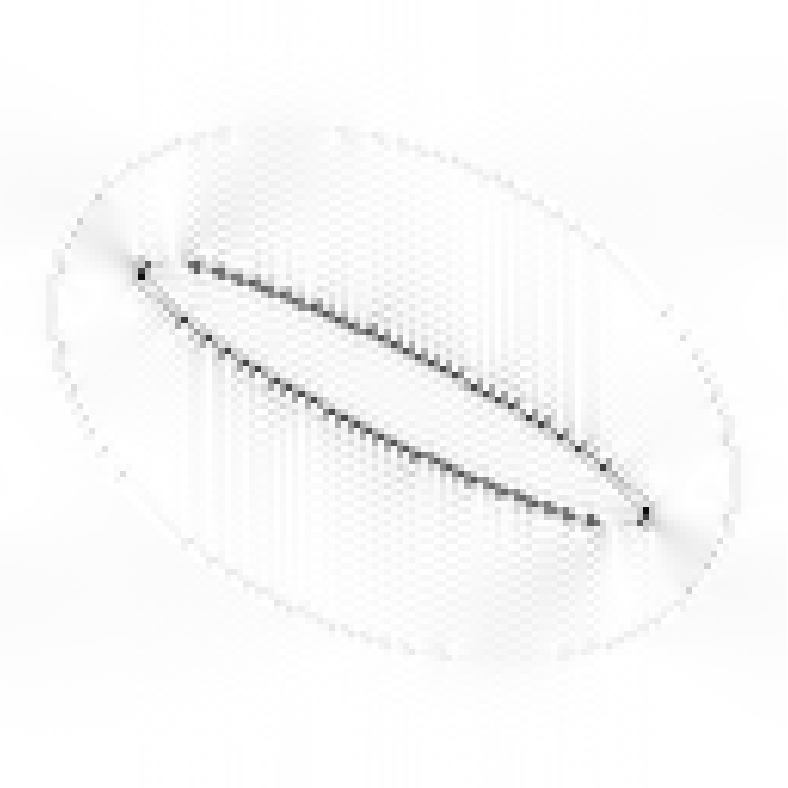}\\[1mm]
    \small
    \begin{tabular}{c l}
      \multirow{3}{*}{$f=f_{\mat{M}_1,\alpha}$,} & $\alpha = (0.45, 0)^\tT$,\\
      & $e_\eff = 0.0017$,\\
      &$e_{\ell^2} = 0.044$.
    \end{tabular}
  \end{subfigure}
  \begin{subfigure}[t]{.32\textwidth}\centering
    \includegraphics[width=0.95\textwidth]{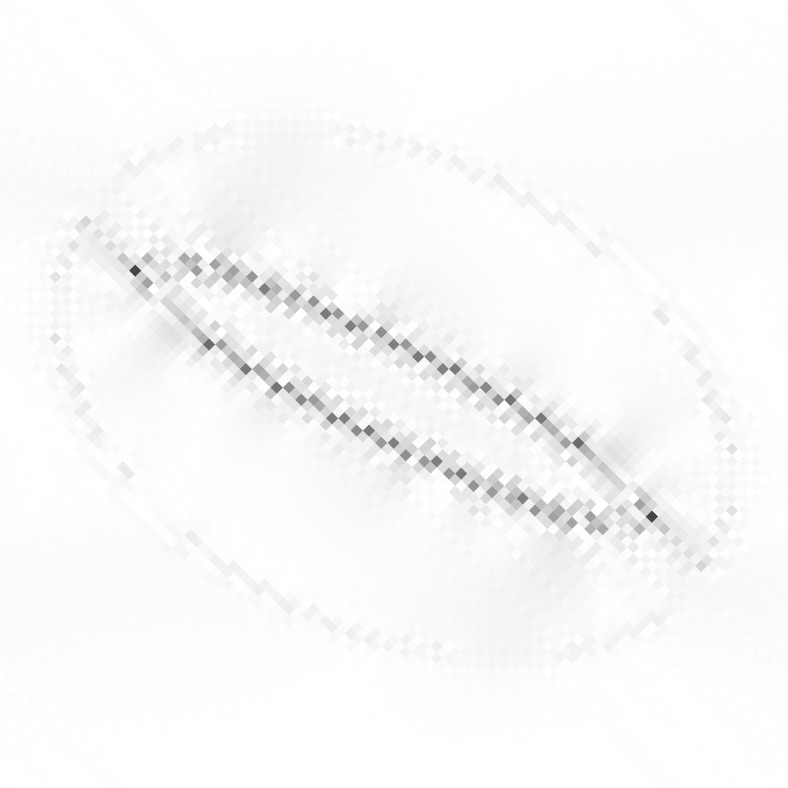}\\[1mm]
    \small
    \begin{tabular}{c l}
      \multirow{3}{*}{$f=f_{\mat{M}_2,\alpha}$,} & $\alpha = (0.25, 0)^\tT$,\\
      & $e_\eff = 0.0013$,\\
      &$e_{\ell^2} = 0.050$.
    \end{tabular}
  \end{subfigure}
  \begin{subfigure}[t]{.32\textwidth}\centering
    \includegraphics[width=0.95\textwidth]{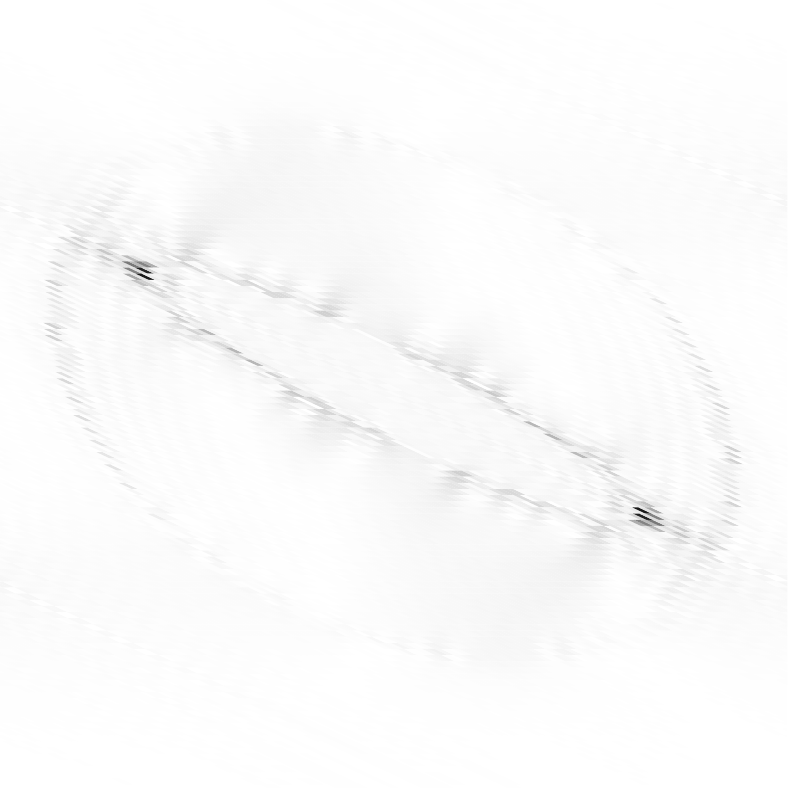}\\[1mm]
    \small
    \begin{tabular}{c l}
      \multirow{3}{*}{$f=f_{\mat{M}_3,\alpha}$,} & $\alpha = (0.4, 0)^\tT$,\\
      & $e_\eff = 0.0024$,\\
      &$e_{\ell^2} = 0.025$.
    \end{tabular}
  \end{subfigure}
    \\[.25\baselineskip]
  \begin{subfigure}[t]{.32\textwidth}\centering
    \includegraphics[width=0.95\textwidth]{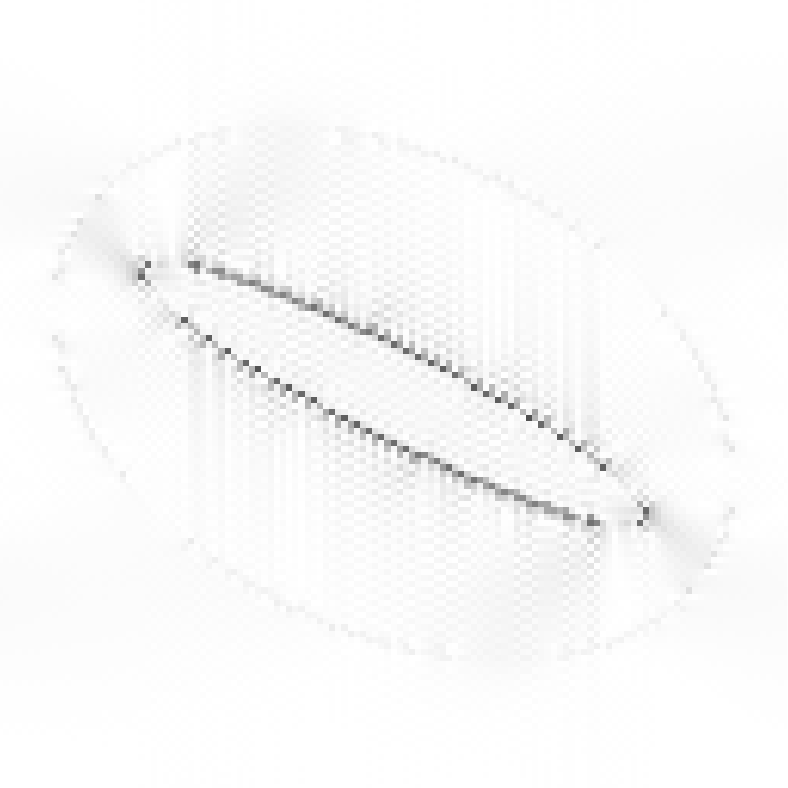}\\[1mm]
    \small
    \begin{tabular}{c l}
      \multirow{2}{*}{$f=D_{\mat{M}_1}$,} & $e_\eff = 0.0038$,\\
      &$e_{\ell^2} = 0.043$.
    \end{tabular}
  \end{subfigure}
  \begin{subfigure}[t]{.32\textwidth}\centering
    \includegraphics[width=0.95\textwidth]{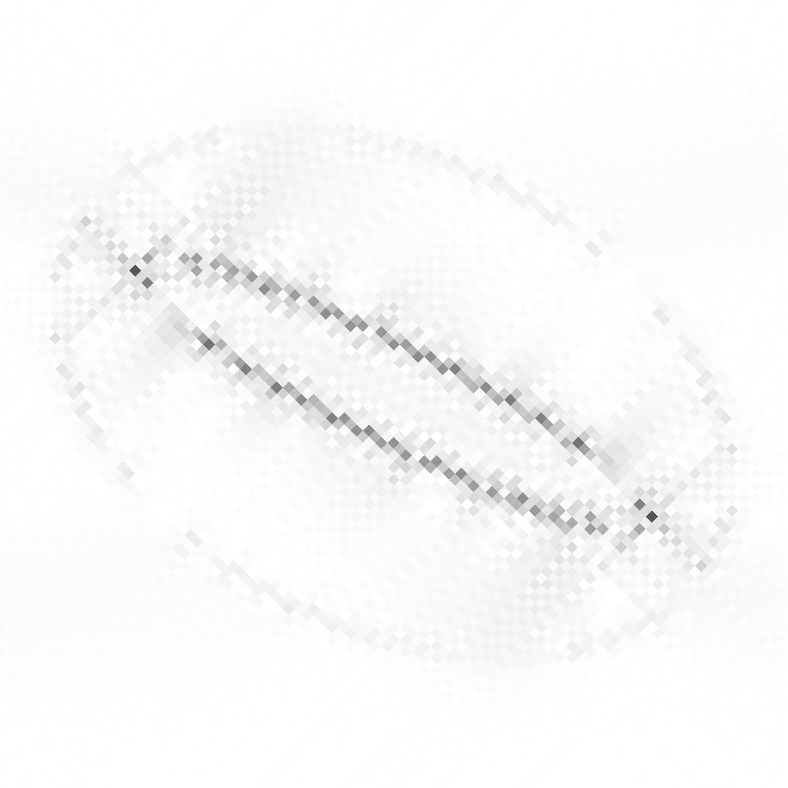}\\[1mm]
    \small
    \begin{tabular}{c l}
      \multirow{2}{*}{$f=D_{\mat{M}_2}$,} & $e_\eff = 0.0034$,\\
      &$e_{\ell^2} = 0.047$.
    \end{tabular}
  \end{subfigure}
  \begin{subfigure}[t]{.32\textwidth}\centering
    \includegraphics[width=0.95\textwidth]{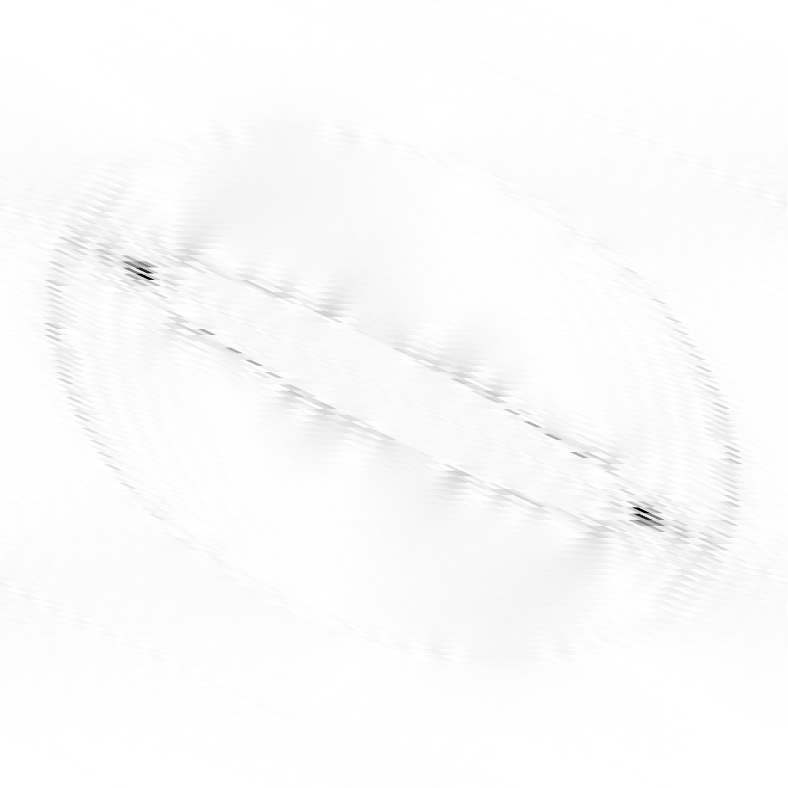}\\[1mm]
    \small
    \begin{tabular}{c l}
      \multirow{2}{*}{$f=D_{\mat{M}_3}$,} & $e_\eff = 0.0036$,\\
      &$e_{\ell^2} = 0.022$.
    \end{tabular}
  \end{subfigure}
  \\[.25\baselineskip]
  \begin{subfigure}[t]{.32\textwidth}\centering
    \includegraphics[width=0.95\textwidth]{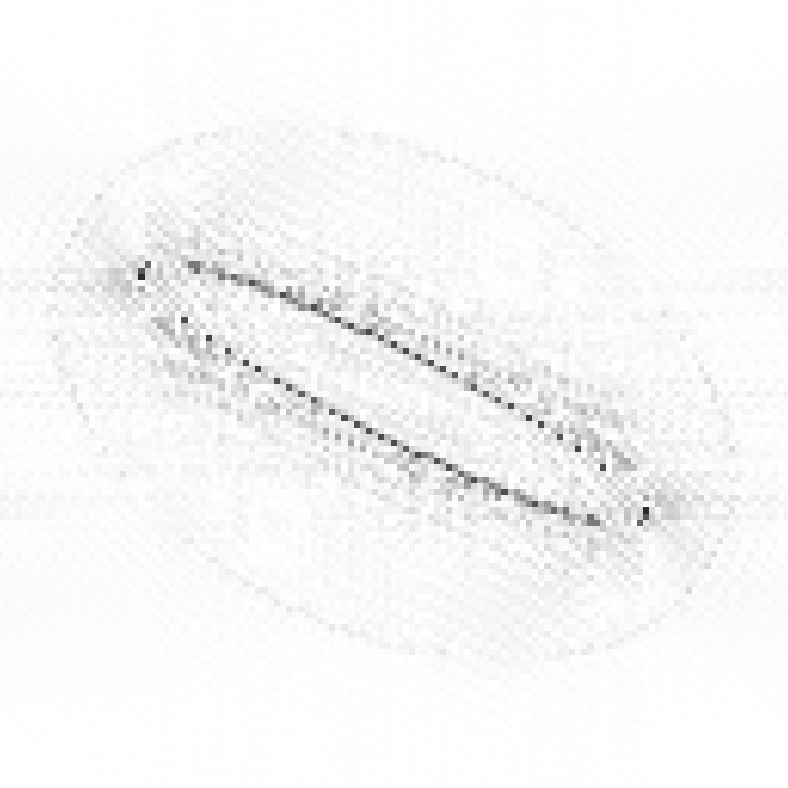}\\[1mm]
    \small
    \begin{tabular}{c l}
      \multirow{2}{*}{$f=B_{\mat{M}_1,\Xi}$,} & $e_\eff = 0.0026$,\\
      &$e_{\ell^2} = 0.052$.
    \end{tabular}
  \end{subfigure}
  \begin{subfigure}[t]{.32\textwidth}\centering
    \includegraphics[width=0.95\textwidth]{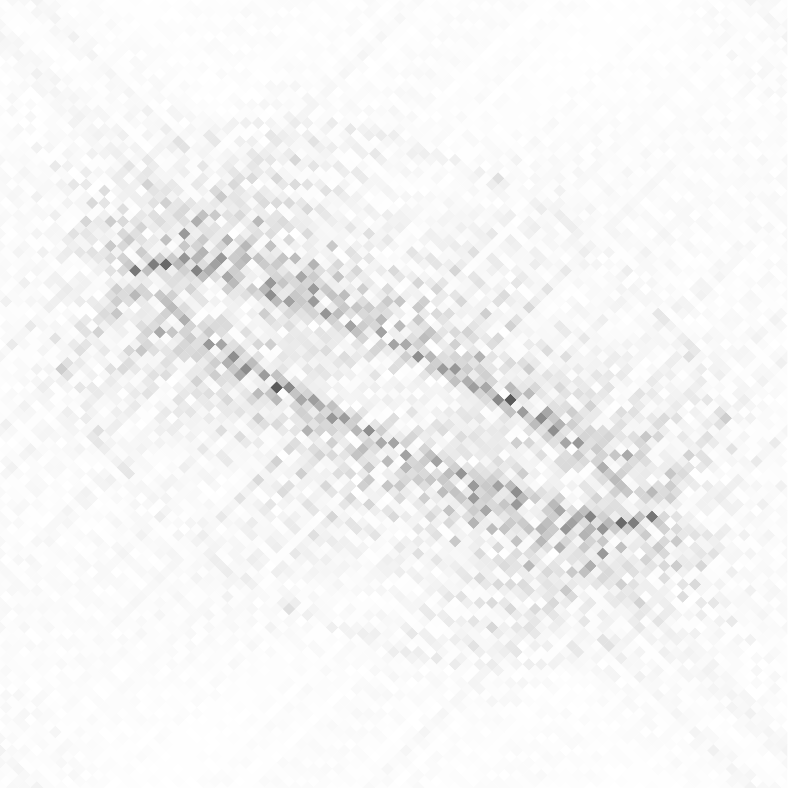}\\[1mm]
    \small
    \begin{tabular}{c l}
      \multirow{2}{*}{$f=B_{\mat{M}_2,\Xi}$,} & $e_\eff = 0.0028$,\\
      &$e_{\ell^2} = 0.059$.
    \end{tabular}
  \end{subfigure}
  \begin{subfigure}[t]{.32\textwidth}\centering
    \includegraphics[width=0.95\textwidth]{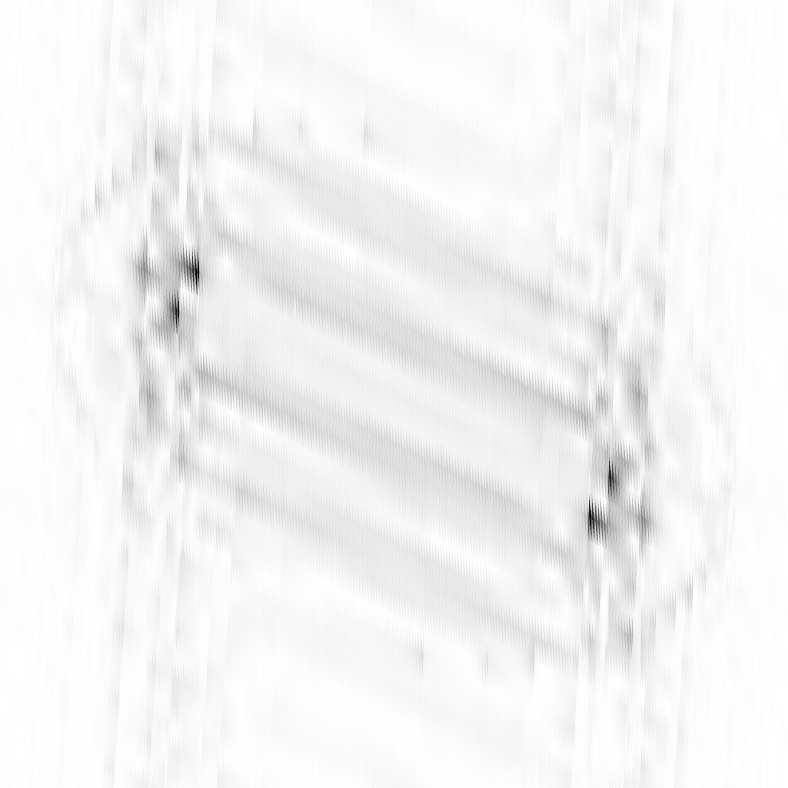}\\[1mm]
    \small
    \begin{tabular}{c l}
      \multirow{2}{*}{$f=B_{\mat{M}_3,\Xi}$,} & $e_\eff = 0.016$,\\
      &$e_{\ell^2} = 0.074$.
    \end{tabular}
  \end{subfigure}
  \caption[]{
    The
    $e_\text{log}$-error of the strain field $\varepsilon_{11}$ given by
    $e_\text{log} = \text{log}(1+|\varepsilon + \tilde{\varepsilon}|)$ using the
    colour bar at the top. In the first row the solution using the de la Vallée
    Poussin Kernel with optimal slopes (regarding the $e_\eff$-error), in the
    second row using the Dirichlet kernel, and in the third row using the Box
    spline
    $\Xi = \bigl(\begin{smallmatrix}
	1 & 0 & 1 & 0 & 1 & 0 & 1 & 0 \\
	0 & 1 & 0 & 1 & 0 & 1 & 0 & 1 
    \end{smallmatrix}\bigr)$.
    From left to right for the pattern matrices 
    $\mat{M}_1 = 
    \bigl(\begin{smallmatrix} 128 & 0\\0 & 128 \end{smallmatrix}\bigl)$,
    $\mat{M}_2 = 
    \bigl(\begin{smallmatrix} 64 & 64\\-64 & 64 \end{smallmatrix}\bigl)$,
    $\mat{M}_3 = 
    \bigl(\begin{smallmatrix} 128 & 272\\0 & 128 \end{smallmatrix}\bigl)$.
    The Box spline plot for $\mat{M}_3$ is an outlier and the scale is adjusted
    by a factor of $1.5$ to account for the large error.
  }\label{fig:Hashin:Ex}
\end{figure}

  Figure~\ref{fig:Hashin:Ex} shows the logarithmic error of the analytic solution
  $\tilde{\varepsilon}$ to solution $\varepsilon$, i.e. $e_\text{log} \coloneqq
  \log(1+\abs{\tilde{\varepsilon}_{11} - \varepsilon_{11}})$. The top row shows
  the error corresponding to the de la Vall\'ee Poussin means with parameters
  $\alpha_1$ and $\alpha_2$ which give the smallest $e_\eff$.  The middle row
  shows the error using the Dirichlet kernel. For illustration purposes each
  pixel has the form of the corresponding unit cell $\mat{M}^{-1}
  [-\frac12,\frac12]^2$ centred at each pattern point $\vect{y} \in
  \Pattern(\mat{M})$.

  The value of the relative
  $\ell_2$-error is given by the formula
  \begin{equation*}
    e_{\ell^2} \coloneqq
    \norm{\varepsilon - \tilde{\varepsilon}} \norm{\tilde{\varepsilon}}^{-1}.
  \end{equation*}

  While the error in the effective stiffness matrix can be drastically improved
  using de la Vall\'ee Poussin means, the $\ell^2$-error gets larger, however,
  only be a few percent. The strain field stemming from the numerical
  computation using de la Vall\'ee Poussin means shows less Gibbs phenomena as
  can be seen most prominently in the solution for pattern matrix $\mat{M}_2$.

  The decrease of $e_\eff$ is caused by the smoothing of the Gibbs phenomena
  around discontinuities of the solution with higher $\alpha_1$.
  With smaller de la Vall\'ee Poussin parameters $\alpha_1$ and $\alpha_2$ the
  polynomial reproduction is better and thus interfaces (edges) are sharper. In
  total they introduce a trade off between damping Gibbs phenomena and
  sharpness of the interfaces.

\subsection{Periodised Box Splines}
Finite elements with one quadrature point are directly included in this framework. They are obtained by choosing a suitable Box spline~\(f=B_{\mat{M},\Xi}\) as ansatz function for
the space of translates~\(V_{\mat{M}}^f\).
In the bottom row of Figure~\ref{fig:Hashin:Ex} the logarithmic error between
the numerical solution and the analytical solution is depicted, using the
matrices $\mat{M}_1$, $\mat{M}_2$, $\mat{M}_3$ from above. For all
computations the Box spline~$B_{\mat{M},\Xi}$ with
\begin{equation}
  \Xi = \begin{pmatrix}
    1 & 0 & 1 & 0 & 1 & 0 & 1 & 0 \\
    0 & 1 & 0 & 1 & 0 & 1 & 0 & 1 
  \end{pmatrix}
\end{equation}
is used. This ansatz function corresponds to a finite element of third order
with reduced integration. The bracket sums in~\eqref{eq:gamma_tilde_fourier} are
precomputed with $33^2$ terms each.

A comparison of the relative error in the effective matrix~$e_\eff$ between the
Box splines (bottom row) and the Dirichlet kernel (middle row) shows that for
matrices $\mat{M}_1$ and $\mat{M}_2$ the error in the effective matrix can be
reduced by about $25\%$, whereas the error is approximately tripled for
$\mat{M}_3$. The $\ell^2$-error $e_{\ell^2}$ is increased in all three cases.
The pattern matrix $\mat{M}_3$  was optimized for the Dirichlet kernel with
respect to the $\ell^2$-norm. For Box splines this pattern leads to an aliasing
effect reducing the quality of the solution.

The choice of the Box spline was not optimized and the quality of the effective
matrix might be further increased by a different choice of $\Xi$. This
corresponds to a fine tuning with respect to dominant directions in the pattern
unit cell $\mat{M}^{-1} \bigl[-\frac12,\frac12\bigr)^d$.

\section{Conclusion}\label{sec:summary}

The introduced framework unifies and analyses for the first time the truncated
Fourier series approach and finite element ansatz functions. In this framework
the periodised Green operator possesses the same properties as the Green
operator of the Galerkin method of Vond\v{r}ejc et.al~\cite{Vondrejc2014}. The
projection operator $\Gamma^0 \Stiffness^0$ emerges as a special case of
Dirichlet kernel translates. The periodised Green operator $\Gamma^\p
\Stiffness^0$ can further be characterised to be a projection if and only if the
space $V_\mat{M}^f$ is the one derived from a Dirichlet kernel.
The finite element method emerges for certain Box splines and thus the constant
finite elements of Brisard and Dormieux~\cite{BrisardDormieux2010Framework} are
included. For these Box splines the infinite Bracket sums have to be precomputed
up to a given precision, but the algorithm has the same complexity as the
Fourier framework. 

Finite elements with more sophisticated quadrature rules can also be viewed in
terms of this framework. However, their performance is an open question for
future work. How to choose a certain Box spline especially with respect to anisotropies present in the data is a point for further analysis. Convergence of
the numerical discretisation towards the continuous case and the convergence of
the algorithm will be dealt with in a paper on our road map. Both seem to be
strongly suggested by the numerical results. Finally the periodic
multiresolution analysis can be used to extend this framework in order to exploit sparsity properties of given data. This is also a point for future work.

\paragraph{Acknowledgement} The authors would like to thank Bernd Simeon and
Gabriele Steidl for their valuable comments on preliminary versions of this manuscript and fruitful discussions.

\printbibliography
\end{document}